\documentclass[11pt]{article}
\usepackage[a4paper,margin=2.5cm]{geometry}
\usepackage[T1]{fontenc}
\usepackage{lmodern}
\usepackage[utf8]{inputenc}
\usepackage{amsmath,amssymb,bm,mathtools,mathrsfs}
\usepackage{booktabs,array,microtype,graphicx}
\usepackage{adjustbox}
\usepackage{algorithm,algpseudocode,placeins}
\usepackage{amsthm,xcolor,url}
\usepackage[font=small,labelfont=bf,skip=5pt]{caption}
\usepackage[colorlinks,linkcolor=blue,citecolor=blue,urlcolor=blue,breaklinks]{hyperref}
\numberwithin{equation}{section}
\newtheorem{proposition}{Proposition}[section]
\newenvironment{keywords}{\par\medskip\noindent\small\textbf{Keywords.}~}{\par}
\newenvironment{MSCcodes}{\par\smallskip\noindent\small\textbf{MSC codes.}~}{\par}
\usepackage{pgfplots,pgfplotstable}
\pgfplotsset{compat=1.18}
\usetikzlibrary{arrows.meta}
\newcommand{\ii}{\mathrm{i}}
\newcommand{\R}{\mathbb{R}}

\newcommand{\vct}[1]{\bm{#1}}
\newcommand{\E}{\bm E}
\newcommand{\f}{\bm f}
\newcommand{\x}{\bm x}
\newcommand{\tht}{\bm\theta}
\newcommand{\dvec}{\bm d}
\newcommand{\svec}{\bm s}
\newcommand{\cc}{\bm c}
\newcommand{\ppw}{\mathrm{ppw}}
\newcommand{\ph}{\mathrm{ph}}
\newcommand{\pml}{\mathrm{pml}}
\newcommand{\curl}{\nabla\!\times}

\newcommand{\diag}{\operatorname{diag}}
\newcommand{\tr}{\operatorname{tr}}

\newcommand{\norm}[1]{\lVert#1\rVert}

\newcommand{\calS}{\mathcal S}

\newcommand{\calV}{\mathcal V}
\newcommand{\calK}{\mathcal K}
\newcommand{\calB}{\mathcal B}
\newcommand{\doi}[1]{\href{https://doi.org/#1}{\nolinkurl{#1}}}
\hypersetup{pdftitle={Toward Extreme-Scale Full-Wave Maxwell Simulation: Low-Dispersion Discretization and a Massively Parallel Three-Grid Solver},pdfauthor={Huangxin Chen, Shubin Fu, Zixiao Zhao}}
\date{}
\title{Toward Extreme-Scale Full-Wave Maxwell Simulation:
Low-Dispersion Discretization and a Massively Parallel Three-Grid Solver}
\author{Huangxin Chen\thanks{School of Mathematical Sciences and Fujian Provincial Key Laboratory on Mathematical Modeling and High Performance Scientific Computing, Xiamen University, Fujian 361005, China (\href{mailto:chx@xmu.edu.cn}{chx@xmu.edu.cn}). Huangxin Chen was supported by National Key Research and Development Project of China (Grant No. 2024YFA1012600), National Natural Science Foundation of China (Grant No. 12471345, 12526102), and Fujian Provincial Natural Science Foundation of China (Grant No. 2026J011002).}
\and Shubin Fu\thanks{Corresponding author. School of Mathematical Sciences, Eastern Institute of Technology, Ningbo, Zhejiang 315200, P. R. China (\href{mailto:sfu@eitech.edu.cn}{sfu@eitech.edu.cn}).}
\and Zixiao Zhao\thanks{School of Mathematical Sciences, Shanghai Jiao Tong University, Shanghai 200240, P. R. China; School of Mathematical Sciences, Eastern Institute of Technology, Ningbo, Zhejiang 315200, P. R. China.}}
\begin{document}
\maketitle

\begin{abstract}
Full-wave simulation of three-dimensional time-harmonic Maxwell problems
over hundreds of wavelengths leads to indefinite systems with tens to
hundreds of billions of unknowns. We present a matrix-free three-grid
solver for such systems and its implementation on GPUs
and CPUs. The solver is built on an interpolated optimized
finite-difference (IOFD) discretization on the Yee grid, whose compact
face and edge filters keep the phase error small at 6 points per
wavelength, and it applies the same compact stencil on all three grids:
an unshifted intermediate-grid correction is preconditioned by a
complex-shifted auxiliary cycle, so that one preconditioner application
performs a fixed number of stencil sweeps and reductions with storage
linear in the number of unknowns. The implementation uses a partition
nested across the grids, $s$-step smoothers that need one global
reduction per smoothing call, and bfloat16 storage of the smoother bases.
On 216 GPUs, complete solves with approximately $71.7$
billion complex unknowns take 77--504 seconds in homogeneous and
heterogeneous media, with weak-scaling efficiencies of 89--91\% per
preconditioner application; strong scaling reduces solves with 2.7
billion unknowns to 4.5--18 seconds on 64 GPUs. A CPU implementation
reaches approximately $574$ billion unknowns on 13\,824 cores, and the
computed fields agree with the outgoing Green field over domains of up
to 956 wavelengths per direction.
\end{abstract}

\begin{keywords}
time-harmonic Maxwell equations, matrix-free three-grid solver,
interpolated optimized finite differences, GPU computing,
high-performance computing, weak and strong scaling
\end{keywords}

\begin{MSCcodes}
65N06, 65F08, 65N55, 65Y05, 78M20
\end{MSCcodes}

\section{Introduction}\label{sec:intro}

Time-harmonic Maxwell equations govern electromagnetic propagation in
scattering, antenna design, photonics, and imaging \cite{Monk2003}. When
the domain spans many wavelengths, two requirements meet. A small error in
the numerical wavenumber accumulates into a large phase error along long
propagation paths, while refinement is expensive: with $s$ points per
wavelength, a domain of $K$ wavelengths in each direction requires
$\mathcal O((sK)^3)$ unknowns. The resulting systems are indefinite and,
with perfectly matched layers (PMLs), complex and anisotropic
\cite{Berenger1994,ChewWeedon1994}. Simulation over hundreds of
wavelengths therefore needs a discretization that is accurate at modest
resolution and a solver whose work and storage grow linearly with the
number of unknowns and map efficiently onto large GPU and CPU
systems.

On the discretization side, the Yee grid \cite{Yee1966} and
N{\'e}d{\'e}lec elements \cite{Nedelec1980} preserve the curl structure of
Maxwell equations but not the phase velocity. Phase accuracy has been
improved through element order and quadrature blending
\cite{Ainsworth2004,AinsworthWajid2010,WajidAyub2012}, compact
high-order staggered schemes \cite{VersanoEtAl2024}, and weights or
stabilization chosen from the discrete dispersion relation
\cite{RewienskiMrozowski2000,JagalurMohan2013,CocquetGander2026}. For the
Helmholtz equation, optimized compact differences \cite{JoShinSuh1996} and
Stolk's interpolated optimized finite differences (IOFD)
\cite{Stolk2016} fit the stencil at the finite resolution that is actually
used and tabulate the coefficients as functions of the local wavenumber.
For Maxwell equations such a design must control two transverse
polarizations while retaining the discrete gradient null space, and, for
large-scale use, it must keep the compact stencil that a matrix-free solver
applies on every level.

On the solver side, time-harmonic Maxwell systems have been solved by
sparse direct methods \cite{AmestoyEtAl2001} and their
block low-rank variants \cite{AmestoyEtAl2015}, by overlapping Schwarz methods
\cite{BonazzoliEtAl2019} and sweeping preconditioners
\cite{TsujiEtAl2012,TsujiYing2012,PoulsonEtAl2013}, which transport information
through coarse problems or approximate factorizations, and by geometric
multigrid \cite{GopalakrishnanEtAl2004,ChanaudEtAl2014}. For multigrid at
high frequency, the coarse grids must represent the propagating waves
while their solution stays inexpensive. Complex-shifted multigrid \cite{ErlanggaEtAl2006,
CoolsVanroose2013}, dispersion-matched coarse operators
\cite{StolkEtAl2014}, and inexact unshifted coarse solves preconditioned
by a shifted cycle \cite{CalandraEtAl2013} address these requirements; the
three-grid methods of \cite{FuEtAlHelmholtz,FuMaxwell} combine them into a
physical and an auxiliary hierarchy for IOFD Helmholtz and edge-element
Maxwell discretizations.

In this paper, we develop an extreme-scale three-grid solver for the
Maxwell IOFD discretization, in which all three levels use the same compact IOFD stencil and are applied
matrix-free, so one preconditioner application performs a fixed number of
stencil sweeps, vector updates, and reductions, with storage linear in the
number of unknowns. Since the coarse levels carry most Krylov steps but
little work, their interface exchanges and reductions decide parallel
efficiency; we organize the implementation around this observation. The
main contributions are as follows.
\begin{itemize}
\item A Maxwell IOFD discretization on the Yee grid: a six-parameter family
of compact face and edge filters that retains the discrete curl and its
gradient null space, with coefficients fitted to both polarizations at
finite resolution and selected from the local material. At 6 points per
wavelength it reduces the maximum wavenumber error of standard Yee by a
factor of about 1700.
\item A matrix-free three-grid solver and its implementation on GPUs
and CPUs, with a partition nested across the grids, $s$-step
smoothers that need one reduction per smoothing call, and bfloat16
storage of the smoother bases, which leaves the convergence and the field
accuracy unchanged.
\item A performance study at extreme scale: complete solves with
approximately $71.7$ billion complex unknowns, spanning 476 wavelengths in
the homogeneous case, in 77--504 seconds on 216 GPUs with
per-application weak-scaling efficiencies of 89--91\%; strong scaling that
reduces billion-unknown solves to seconds; and CPU solves with about
$574$ billion unknowns on 13\,824 cores, validated against the outgoing Green
field over 956 wavelengths.
\end{itemize}

The rest of this paper is organized as follows. In
Section~\ref{sec:discretization}, we construct the Maxwell IOFD
discretization and optimize its coefficients for both polarizations. On
this discretization we build the three-grid preconditioner and its
distributed implementation in Section~\ref{sec:threegrid}. In
Section~\ref{sec:field-validation}, we validate the computed wavefields
against the outgoing Green field and select the stopping tolerance, and
in Section~\ref{sec:gpu-results} we report weak and strong scaling on
GPUs and CPUs. The dispersion analysis and the proofs are given in the
appendices.

\section{Maxwell IOFD Discretization}\label{sec:discretization}

In this section, we
construct the Maxwell discretization that the solver of
Section~\ref{sec:threegrid} applies matrix-free. We define compact face
and edge filters on the Yee grid (Section~\ref{sec:construction}),
optimize their coefficients for both polarizations at the operating
resolution (Section~\ref{sec:symbol}), add a source and field filter
(Section~\ref{sec:amplitude}), and select the coefficients locally in
variable media with a constitutive PML (Section~\ref{sec:pml}). Derivations, the low-wavenumber expansion,
and the error bounds behind the coefficient design are collected in
Appendix~\ref{sec:dispersion}.

\subsection{The Maxwell problem}\label{sec:model}

With time dependence $\exp(-\ii\omega t)$, the electric field generated by
an impressed source satisfies
\begin{equation}
 \curl\bigl(\mu(\x)^{-1}\curl\E\bigr)
       -k_0^2\epsilon(\x)\E=\f
       \quad\text{in }\R^3,
 \label{eq:variable-pde}
\end{equation}
where $\epsilon$ and $\mu$ are positive relative material parameters and
$k_0$ is the reference wavenumber; the source is supported in a bounded
physical region $\Omega_{\ph}$, and the field is outgoing at infinity
\cite{Monk2003}. For constant material values the equation reduces to
\begin{equation}
 \curl\curl\E-k^2\E=\f_0,
 \qquad k=k_0\sqrt{\epsilon\mu},
 \label{eq:pde}
\end{equation}
with $\f_0=\mu\f$, whose two transverse plane-wave polarizations share the
dispersion sphere $|\vct\xi|=k$. The discretization should reproduce this
surface for both polarizations at a prescribed resolution, measured by
\begin{equation}
 \eta=\frac{kh}{2\pi}=\frac1{\ppw},\qquad \kappa=kh,
 \label{eq:kappa}
\end{equation}
where $h$ is the mesh spacing and $\ppw$ the number of points per wavelength.

\subsection{A compact face--edge construction}\label{sec:construction}

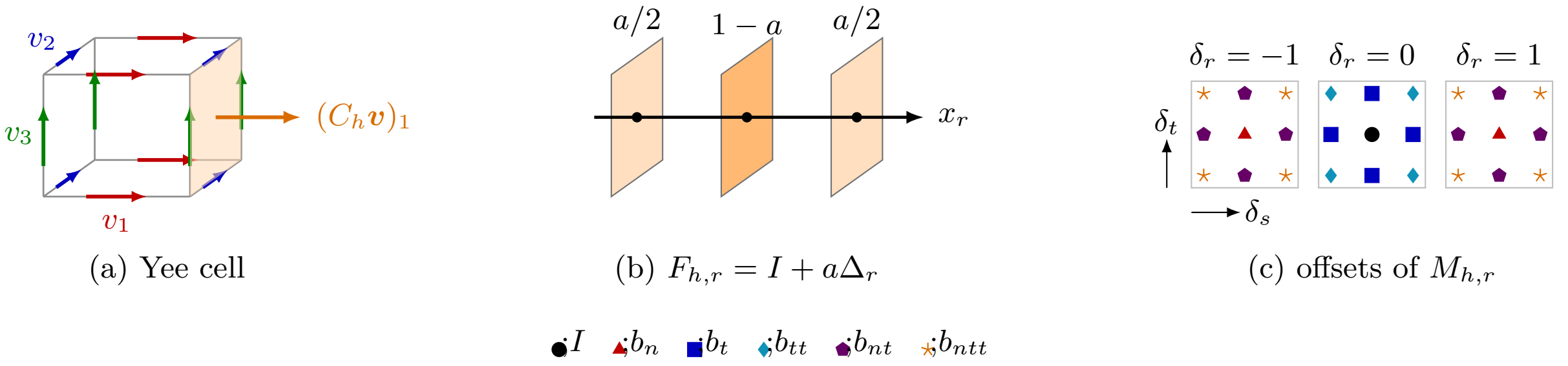
\begin{figure}[htbp]
\centering
\tikzset{
 obl/.style={x={(1cm,0cm)},y={(0.42cm,0.30cm)},z={(0cm,1cm)}},
 cube/.style={black!45,line width=0.5pt},
 ex/.style={-{Latex[length=1.6mm]},line width=1.1pt,red!75!black},
 ey/.style={-{Latex[length=1.6mm]},line width=1.1pt,blue!75!black},
 ez/.style={-{Latex[length=1.6mm]},line width=1.1pt,green!50!black},
}
\begin{tikzpicture}[obl,scale=1.25,font=\footnotesize]
 \foreach \y in {0,1} \foreach \z in {0,1} \draw[cube] (0,\y,\z)--(1.2,\y,\z);
 \foreach \x in {0,1.2} \foreach \z in {0,1} \draw[cube] (\x,0,\z)--(\x,1,\z);
 \foreach \x in {0,1.2} \foreach \y in {0,1} \draw[cube] (\x,\y,0)--(\x,\y,1);
 \foreach \y in {0,1} \foreach \z in {0,1} \draw[ex] (0.35,\y,\z)--(0.85,\y,\z);
 \foreach \x in {0,1.2} \foreach \z in {0,1} \draw[ey] (\x,0.25,\z)--(\x,0.75,\z);
 \foreach \x in {0,1.2} \foreach \y in {0,1} \draw[ez] (\x,\y,0.25)--(\x,\y,0.75);
 \fill[orange!30,opacity=0.55] (1.2,0,0)--(1.2,1,0)--(1.2,1,1)--(1.2,0,1)--cycle;
 \draw[-{Latex[length=1.8mm]},line width=1.1pt,orange!85!black] (1.2,0.5,0.5)--(1.9,0.5,0.5)
   node[right]{$(C_h\vct v)_1$};
 \node[red!75!black] at (0.6,0,-0.22) {$v_1$};
 \node[blue!75!black] at (-0.22,0.5,1.12) {$v_2$};
 \node[green!50!black] at (-0.2,0,0.5) {$v_3$};
 \node at (0.8,0.5,-0.75) {(a) Yee cell};
\end{tikzpicture}\hfill
\begin{tikzpicture}[obl,scale=1.25,font=\footnotesize]
 \foreach \x/\w/\c in {0/{a/2}/orange!25,0.9/{1-a}/orange!55,1.8/{a/2}/orange!25}{
   \fill[\c] (\x,0,0)--(\x,1,0)--(\x,1,1)--(\x,0,1)--cycle;
   \draw[black!55] (\x,0,0)--(\x,1,0)--(\x,1,1)--(\x,0,1)--cycle;
   \fill (\x,0.5,0.5) circle (1.2pt);
   \node[above] at (\x,0.5,1.05) {$\w$};}
 \draw[-{Latex[length=1.8mm]},line width=0.9pt] (-0.35,0.5,0.5)--(2.35,0.5,0.5) node[right]{$x_r$};
 \node at (0.9,0.5,-0.75) {(b) $F_{h,r}=I+a\Delta_r$};
\end{tikzpicture}\hfill
\begin{tikzpicture}[scale=0.42,font=\footnotesize]
 \foreach \dn/\X in {-1/0,0/3.1,1/6.2}{
   \draw[black!25] (\X-1.3,-1.3) rectangle (\X+1.3,1.3);
   \foreach \ds in {-1,0,1} \foreach \dt in {-1,0,1}{
     \pgfmathtruncatemacro\nt{abs(\ds)+abs(\dt)}\pgfmathtruncatemacro\nn{abs(\dn)}
     \ifnum\nn=0 \ifnum\nt=0 \def\col{black}\def\mk{*}\fi
                 \ifnum\nt=1 \def\col{blue!75!black}\def\mk{square*}\fi
                 \ifnum\nt=2 \def\col{cyan!70!black}\def\mk{diamond*}\fi
     \else       \ifnum\nt=0 \def\col{red!75!black}\def\mk{triangle*}\fi
                 \ifnum\nt=1 \def\col{violet!80!black}\def\mk{pentagon*}\fi
                 \ifnum\nt=2 \def\col{orange!85!black}\def\mk{star}\fi
     \fi
     \node[\col,inner sep=0pt] at (\X+\ds,\dt) {\pgfuseplotmark{\mk}};}
   \node[above] at (\X,1.3) {$\delta_r=\dn$};}
 \draw[-{Latex[length=1.4mm]}] (-1.3,-1.9)--(-0.1,-1.9) node[right,inner sep=1pt]{$\delta_s$};
 \draw[-{Latex[length=1.4mm]}] (-1.9,-1.3)--(-1.9,-0.1) node[above,inner sep=1pt]{$\delta_t$};
 \node at (3.1,-3.3) {(c) offsets of $M_{h,r}$};
\end{tikzpicture}

\vspace{2pt}
{\footnotesize
\tikz{\node[black,inner sep=0] at (0,0){\pgfuseplotmark{*}};};$I$\quad
\tikz{\node[red!75!black,inner sep=0] at (0,0){\pgfuseplotmark{triangle*}};};$b_n$\quad
\tikz{\node[blue!75!black,inner sep=0] at (0,0){\pgfuseplotmark{square*}};};$b_t$\quad
\tikz{\node[cyan!70!black,inner sep=0] at (0,0){\pgfuseplotmark{diamond*}};};$b_{tt}$\quad
\tikz{\node[violet!80!black,inner sep=0] at (0,0){\pgfuseplotmark{pentagon*}};};$b_{nt}$\quad
\tikz{\node[orange!85!black,inner sep=0] at (0,0){\pgfuseplotmark{star}};};$b_{ntt}$}
\caption{Unknowns and filter supports of the face--edge construction. (a) The edge unknowns $v_1,v_2,v_3$ of one Yee cell and a face value $(C_h\vct v)_1$, formed from the four edges around its face. (b) The face filter acts along the face normal with weights $a/2$, $1-a$, $a/2$. (c) The edge filter of component $r$ couples the $27$ parallel edges at offsets $(\delta_r,\delta_s,\delta_t)\in\{-1,0,1\}^3$, drawn as three layers $\delta_r=-1,0,1$ along the edge direction, with every edge pointing along $x_r$; the six symbols mark the symmetry classes, which are first reached by the terms $I$, $b_n\Delta_r$, $b_t(\Delta_s+\Delta_t)$, $b_{tt}\Delta_s\Delta_t$, $b_{nt}\Delta_r(\Delta_s+\Delta_t)$, and $b_{ntt}\Delta_r\Delta_s\Delta_t$ of \eqref{eq:edge}.}
\label{fig:stencil}
\end{figure}

The $r$th electric component is sampled at
$h(\vct j+\tfrac12\vct e_r)$, $\vct j\in\mathbb Z^3$, and the discrete
curl $C_h$ maps edge fields to face fields, for example
\begin{equation}
 (C_h\vct v)_{1,i,j+1/2,l+1/2}
 =\frac{v_{3,i,j+1,l+1/2}-v_{3,i,j,l+1/2}}{h}
  -\frac{v_{2,i,j+1/2,l+1}-v_{2,i,j+1/2,l}}{h},
 \label{eq:curlstencil}
\end{equation}
with the other components by cyclic permutation; together with the
nodal-to-edge gradient $G_h$ these differences satisfy $C_hG_h=0$. The
filters act on one component at a time in the index coordinates of its
own lattice. With the translation $\mathcal T_j$ by one mesh interval in
direction $j$ and the second difference
\begin{equation}
 \Delta_j=\tfrac12(\mathcal T_j+\mathcal T_j^{-1})-I,
 \label{eq:delta}
\end{equation}
the face filter acts along the face normal and the edge filter
componentwise,
\begin{equation}
 F_{h,r}=I+a\Delta_r,
 \label{eq:face}
\end{equation}
\begin{align}
 M_{h,r}={}&I+b_n\Delta_r+b_t(\Delta_s+\Delta_t)
       +b_{tt}\Delta_s\Delta_t \notag\\
       &+b_{nt}\Delta_r(\Delta_s+\Delta_t)
       +b_{ntt}\Delta_r\Delta_s\Delta_t,
       \qquad \{r,s,t\}=\{1,2,3\}.
 \label{eq:edge}
\end{align}
The edge filter is the general $3\times3\times3$ stencil that is even in
each coordinate, treats the two transverse directions equally, and
reproduces constants. With
\begin{equation}
 \cc=(a,b_n,b_t,b_{tt},b_{nt},b_{ntt}),\qquad
 F_h=\diag_rF_{h,r},\quad M_h=\diag_rM_{h,r},
 \label{eq:parameters}
\end{equation}
the homogeneous operator is
\begin{equation}
 A_h=C_h^*F_hC_h-k^2M_h.
 \label{eq:operator}
\end{equation}
Standard Yee is the member $\cc=0$. Each component has a 27-point mass
stencil, and the curl term couples the three edge lattices within one
cell, so the six coefficients shape the dispersion relation without
enlarging the stencil that the solver applies.
Figure~\ref{fig:stencil} shows the unknowns of one Yee cell and the
supports of the two filters.

\begin{proposition}[Consistency and discrete gradient constraint]
\label{prop:consistency}
For constant $k>0$ and bounded coefficients $\cc$, the interior operator
\eqref{eq:operator} is second-order consistent with \eqref{eq:pde}.
Its stiffness $K_h=C_h^*F_hC_h$ satisfies $K_hG_h=0$ and $G_h^*K_h=0$.
Consequently, a solution of $A_h\vct v_h=\vct b_h$ obeys
\begin{equation}
 -k^2G_h^*M_h\vct v_h=G_h^*\vct b_h.
 \label{eq:discretegauss}
\end{equation}
\end{proposition}
The proof is given in Appendix~\ref{sec:proofs}.

\subsection{Dispersion relation and coefficient optimization}\label{sec:symbol}

On a staggered Fourier mode
\begin{equation}
 (v_h)_r(\vct j)=p_r
 \exp\!\bigl(\ii\tht\cdot(\vct j+\tfrac12\vct e_r)\bigr),
 \qquad\tht\in[-\pi,\pi]^3,
 \label{eq:fourier}
\end{equation}
let $s_j=2\sin(\theta_j/2)$, $z_j=\cos\theta_j-1$, and
$S=[\svec]_\times$. The curl symbol is $\ii S/h$, and the filter symbols
$\widehat F=\diag(f_r)$ and $\widehat M=\diag(m_r)$ are
\begin{align}
 f_r={}&1+az_r,\label{eq:Fsymbol}\\
 m_r={}&1+b_nz_r+b_t(z_s+z_t)+b_{tt}z_sz_t\notag\\
       &+b_{nt}z_r(z_s+z_t)+b_{ntt}z_rz_sz_t,
 \label{eq:Msymbol}
\end{align}
so that the dimensionless symbol is
\begin{equation}
 \mathscr P(\tht;\eta)=h^2\widehat A_h(\tht;\eta)
       =\widehat K-\kappa^2\widehat M,
 \qquad\widehat K=S^T\widehat F S.
 \label{eq:dimensionlesssymbol}
\end{equation}

\begin{proposition}[Admissibility and transverse dispersion branches]
\label{prop:branches}
Fix a real coefficient vector $\cc$, and let $f_r$ and $m_r$ be its filter
symbols in \eqref{eq:Fsymbol}--\eqref{eq:Msymbol}. Suppose that
\begin{equation}
 f_r(\tht)>0,\qquad m_r(\tht)>0,
 \quad r=1,2,3,\quad\tht\in[-\pi,\pi]^3.
 \label{eq:positivity}
\end{equation}
For $\tht\ne0$, the generalized eigenvalues of
$(\widehat K,\widehat M)$ are $0,\lambda_1,\lambda_2$, with
$0<\lambda_1\le\lambda_2$. The zero eigenvector is $\svec$, and
each propagating eigenvector satisfies $\svec^T\widehat M\vct p=0$.
For $B=\widehat M^{-1/2}\widehat K\widehat M^{-1/2}$,
\begin{equation}
 \lambda_{1,2}
 =\frac{\tr B\mp\sqrt{2\tr(B^2)-(\tr B)^2}}2,
 \label{eq:explicitbranches}
\end{equation}
and
\begin{equation}
 \det\mathscr P=-\kappa^2\det\widehat M\,
     (\lambda_1-\kappa^2)(\lambda_2-\kappa^2).
 \label{eq:characteristicfactorization}
\end{equation}
Condition~\eqref{eq:positivity} is equivalent to positivity of the face
symbols at $z_r\in\{-2,0\}$ and the edge symbols at
$(z_1,z_2,z_3)\in\{-2,0\}^3$.
\end{proposition}
The proof is given in Appendix~\ref{sec:proofs}.

We call coefficients satisfying \eqref{eq:positivity} \emph{admissible}.
Their propagating waves solve $\lambda_b(\tht;\cc(\eta))=\kappa^2$,
$b=1,2$. Along a direction $\dvec\in\mathbb S^2$ we take the
near-physical root $q_b$ of $\lambda_b(q_b\kappa\dvec;\cc(\eta))=\kappa^2$,
so the numerical wavenumber is $q_bk$, and measure the relative
wavenumber error \cite[Section~2.1]{Stolk2016}
\begin{equation}
 \delta_b(\dvec;\eta)=q_b-1,\qquad
 \mathcal E_{\ph}(\eta)
   =\sup_{\dvec\in\mathbb S^2}\max_{b=1,2}|\delta_b(\dvec;\eta)|.
 \label{eq:metrics}
\end{equation}
Over a distance $L$ the plane-wave phases then differ by
\begin{equation}
 \Delta\phi_b=kL\delta_b
       =2\pi\frac{L}{\lambda_w}\delta_b,
 \qquad\lambda_w=\frac{2\pi}k,
 \label{eq:accumulation}
\end{equation}
so a small relative wavenumber error is essential over many wavelengths,
and $\mathcal E_{\ph}$ also bounds the difference of the two polarization
wavenumbers by $2\mathcal E_{\ph}$.

A low-wavenumber expansion (Proposition~\ref{prop:expansion}) shows that
the $O(\kappa^2)$ error of both polarizations cancels in every direction
if and only if $a=b_t=1/6$ and $b_n=0$, whereas $b_{tt}$, $b_{nt}$, and
$b_{ntt}$ enter only at higher order. At practical resolutions these
higher-order terms are not negligible, so we instead fit all six
coefficients at each tabulated $\eta$, extending the phase fitting of
\cite[Section~5.1]{Stolk2016} to both polarizations:
\begin{equation}
 J_{\ph}(\cc;\eta)
    =\sum_{\dvec\in\mathcal D}w_{\dvec}
                   \sum_{b=1}^2|q_b(\dvec;\eta,\cc)-1|^2,
 \qquad w_{\dvec}>0,
 \label{eq:phaseobjective}
\end{equation}
over a weighted set of directions $\mathcal D\subset\mathbb S^2$, with
admissibility imposed through positive lower bounds on the vertex values in
Proposition~\ref{prop:branches}. Between tabulation points $\eta_j$ the
coefficients are interpolated linearly,
\begin{equation}
 \cc(\eta)=(1-t)\cc_j+t\cc_{j+1},\qquad
 t=\frac{\eta-\eta_j}{\eta_{j+1}-\eta_j},
 \label{eq:linearinterpolation}
\end{equation}
which preserves admissibility because the filter symbols are affine in
$\cc$. Propositions~\ref{prop:radial} and \ref{prop:perturbation} bound
the resulting wavenumber error by the symbol defect divided by the radial
slope of the branch, so the interpolation contributes a second-order term
in the table spacing.

The coefficients are tabulated at $\eta_j=0.1,0.125,\dots,0.4$ together
with $1/12$, $1/6$, and $1/3$; outside $[1/12,0.4]$ the filter
coefficients are frozen at the nearest endpoint while the mass term
retains the actual wavenumber.
Table~\ref{tab:dispersion} compares the maximum relative wavenumber error
over both branches, sampled on a $257\times257$ angular grid in the first
octant, which covers all directions by cubic symmetry. At 6 points per wavelength the
optimized filters reduce the error from $5.25\times10^{-2}$ for standard
Yee to $3.08\times10^{-5}$, so by \eqref{eq:accumulation} the accumulated
phase error shrinks by a factor of about 1700; at 3 points per wavelength
IOFD retains near-physical roots in all directions, whereas the axial Yee
dispersion relation has no real root.

\begin{table}[htbp]
\centering
\caption{Maximum relative wavenumber error over directions and both polarizations for standard Yee and Maxwell IOFD. A dash indicates the absence of a real Yee propagation root in some directions.}
\label{tab:dispersion}
\begin{tabular}{@{}ccc@{}}
\toprule
$\ppw$ & Standard Yee & Maxwell IOFD\\
\midrule
8 & $2.7669\!\times\!10^{-2}$ & $3.0602\!\times\!10^{-5}$\\
6 & $5.2465\!\times\!10^{-2}$ & $3.0848\!\times\!10^{-5}$\\
4 & $1.5017\!\times\!10^{-1}$ & $1.2554\!\times\!10^{-4}$\\
3 & --- & $4.5864\!\times\!10^{-4}$\\
\bottomrule
\end{tabular}
\end{table}

Figure~\ref{fig:phasecomparison} plots the maxima of
Table~\ref{tab:dispersion}, and Figure~\ref{fig:angular6} shows the signed
errors of the two polarizations along an angular cut at 6 points per
wavelength; both remain below $3\times10^{-5}$ in magnitude, and the
sampled maxima give $|k_{h,1}-k_{h,2}|/k\le6.2\times10^{-5}$.

\begin{figure}[htbp]
\centering
\begin{minipage}[t]{0.485\linewidth}
\centering
\begin{tikzpicture}
\begin{semilogyaxis}[
 width=\linewidth,height=5.0cm,
 xlabel={$\eta=1/\ppw$},
 ylabel={Max.\ relative wavenumber error},
 xmin=.105,xmax=.35,ymin=1e-5,ymax=.3,
 xtick={.125,.166666666667,.25,.333333333333},
 xticklabels={$1/8$,$1/6$,$1/4$,$1/3$},
 grid=major,
 tick label style={font=\scriptsize},label style={font=\footnotesize},
 legend style={at={(0.5,1.03)},anchor=south,legend columns=2,font=\scriptsize,cells={anchor=west},/tikz/every even column/.append style={column sep=3mm},draw=none},
]
\addplot+[only marks, mark=square*, mark size=2pt] coordinates {
(0.125,0.0276688327845) (0.166666666667,0.0524653776544) (0.25,0.150167078134)
};
\addlegendentry{Standard Yee}
\addplot+[only marks, mark=*, mark size=2pt] coordinates {
(0.125,3.06018294906e-05) (0.166666666667,3.08477590636e-05) (0.25,0.000125539586833) (0.333333333333,0.000458638938613)
};
\addlegendentry{Maxwell IOFD}
\end{semilogyaxis}
\end{tikzpicture}
\caption{Maximum sampled relative wavenumber error at 8, 6, 4, and 3 points per wavelength. The standard-Yee root is absent in some directions at 3 points per wavelength.}
\label{fig:phasecomparison}
\end{minipage}\hfill
\begin{minipage}[t]{0.485\linewidth}
\centering
\begin{tikzpicture}
\begin{axis}[
 width=\linewidth,height=5.0cm,
 xlabel={Polar angle $\alpha$ (degrees)},
 ylabel={$10^{5}(q_b-1)$},
 xmin=0,xmax=90,xtick={0,15,30,45,60,75,90},
 grid=major,
 tick label style={font=\scriptsize},label style={font=\footnotesize},
 legend style={at={(0.5,1.03)},anchor=south,legend columns=2,font=\scriptsize,cells={anchor=west},/tikz/every even column/.append style={column sep=3mm},draw=none},
]
\addplot+[no marks, thick, solid] coordinates {
(0,-2.898704145) (0.25,-2.897228308) (0.5,-2.892802452) (0.75,-2.885431555) (1,-2.875123895)
(1.25,-2.861891059) (1.5,-2.845747912) (1.75,-2.826712593) (2,-2.804806484) (2.25,-2.780054194)
(2.5,-2.75248352) (2.75,-2.722125423) (3,-2.689013986) (3.25,-2.653186377) (3.5,-2.614682805)
(3.75,-2.573546471) (4,-2.529823522) (4.25,-2.483562988) (4.5,-2.434816734) (4.75,-2.383639394)
(5,-2.330088309) (5.25,-2.274223456) (5.5,-2.216107385) (5.75,-2.155805135) (6,-2.09338417)
(6.25,-2.028914291) (6.5,-1.962467556) (6.75,-1.894118201) (7,-1.823942544) (7.25,-1.752018905)
(7.5,-1.678427506) (7.75,-1.603250385) (8,-1.526571293) (8.25,-1.4484756) (8.5,-1.369050197)
(8.75,-1.288383389) (9,-1.206564795) (9.25,-1.123685246) (9.5,-1.039836674) (9.75,-0.9551120064)
(10,-0.8696050605) (10.25,-0.7834104287) (10.5,-0.6966233724) (10.75,-0.6093397104) (11,-0.5216557037)
(11.25,-0.4336679484) (11.5,-0.3454732604) (11.75,-0.2571685587) (12,-0.1688507626) (12.25,-0.08061666621)
(12.5,0.007437165017) (12.75,0.0952145121) (13,0.182619612) (13.25,0.2695572737) (13.5,0.3559329827)
(13.75,0.4416530186) (14,0.5266245631) (14.25,0.6107558073) (14.5,0.6939560597) (14.75,0.776135858)
(15,0.8572070675) (15.25,0.9370829929) (15.5,1.015678478) (15.75,1.092910009) (16,1.168695814)
(16.25,1.242955963) (16.5,1.315612469) (16.75,1.386589379) (17,1.455812872) (17.25,1.523211354)
(17.5,1.588715547) (17.75,1.652258582) (18,1.713776089) (18.25,1.773206284) (18.5,1.830490061)
(18.75,1.885571073) (19,1.938395824) (19.25,1.988913755) (19.5,2.037077329) (19.75,2.08284212)
(20,2.126166902) (20.25,2.167013748) (20.5,2.205348107) (20.75,2.241138923) (21,2.274358725)
(21.25,2.304983745) (21.5,2.33299403) (21.75,2.358373579) (22,2.381110484) (22.25,2.401197086)
(22.5,2.418630157) (22.75,2.433411104) (23,2.445546208) (23.25,2.455046883) (23.5,2.46193001)
(23.75,2.466218299) (24,2.467940732) (24.25,2.467133091) (24.5,2.463838583) (24.75,2.45810859)
(25,2.45000358) (25.25,2.439594187) (25.5,2.426962536) (25.75,2.412203815) (26,2.395428194)
(26.25,2.376763096) (26.5,2.356355911) (26.75,2.334377161) (27,2.311024151) (27.25,2.28652501)
(27.5,2.261142975) (27.75,2.235180517) (28,2.208982642) (28.25,2.1829383) (28.5,2.157478331)
(28.75,2.133067996) (29,2.110192071) (29.25,2.089331213) (29.5,2.070930175) (29.75,2.055361422)
(30,2.0428908) (30.25,2.033653392) (30.5,2.027646027) (30.75,2.024738245) (31,2.024697999)
(31.25,2.027224832) (31.5,2.031982967) (31.75,2.038629008) (32,2.046831915) (32.25,2.056285376)
(32.5,2.066714004) (32.75,2.077875152) (33,2.089557963) (33.25,2.10158091) (33.5,2.113788663)
(33.75,2.126048808) (34,2.138248711) (34.25,2.150292674) (34.5,2.162099439) (34.75,2.173600039)
(35,2.184735967) (35.25,2.195457635) (35.5,2.205723086) (35.75,2.215496911) (36,2.224749339)
(36.25,2.233455493) (36.5,2.241594746) (36.75,2.249150196) (37,2.256108217) (37.25,2.262458084)
(37.5,2.26819165) (37.75,2.27330308) (38,2.277788619) (38.25,2.28164639) (38.5,2.284876235)
(38.75,2.287479554) (39,2.289459192) (39.25,2.290819321) (39.5,2.291565344) (39.75,2.291703815)
(40,2.291242359) (40.25,2.290189605) (40.5,2.288555135) (40.75,2.286349422) (41,2.283583784)
(41.25,2.280270343) (41.5,2.276421986) (41.75,2.272052323) (42,2.267175659) (42.25,2.26180696)
(42.5,2.255961821) (42.75,2.249656447) (43,2.242907614) (43.25,2.235732656) (43.5,2.228149436)
(43.75,2.22017632) (44,2.211832162) (44.25,2.203136274) (44.5,2.194108412) (44.75,2.184768753)
(45,2.175137873) (45.25,2.165236724) (45.5,2.155086623) (45.75,2.144709223) (46,2.134126499)
(46.25,2.123360721) (46.5,2.112434442) (46.75,2.101370473) (47,2.090191861) (47.25,2.078921872)
(47.5,2.067583966) (47.75,2.056201777) (48,2.044799093) (48.25,2.033399828) (48.5,2.022028003)
(48.75,2.010707721) (49,1.99946314) (49.25,1.988318451) (49.5,1.977297848) (49.75,1.966425502)
(50,1.955725533) (50.25,1.945221978) (50.5,1.934938761) (50.75,1.924899662) (51,1.915128278)
(51.25,1.905647993) (51.5,1.896481935) (51.75,1.887652939) (52,1.879183506) (52.25,1.871095757)
(52.5,1.863411387) (52.75,1.856151616) (53,1.84933714) (53.25,1.842988075) (53.5,1.837123899)
(53.75,1.831763392) (54,1.826924577) (54.25,1.822624649) (54.5,1.818879909) (54.75,1.815705689)
(55,1.813116282) (55.25,1.811124859) (55.5,1.80974339) (55.75,1.808982562) (56,1.808851692)
(56.25,1.809358641) (56.5,1.810509727) (56.75,1.812309633) (57,1.814761321) (57.25,1.817865943)
(57.5,1.821622755) (57.75,1.826029031) (58,1.831079983) (58.25,1.836768694) (58.5,1.843086033)
(58.75,1.850020605) (59,1.857558695) (59.25,1.865684221) (59.5,1.874378702) (59.75,1.883621244)
(60,1.893388523) (60.25,1.9036548) (60.5,1.914391945) (60.75,1.925569468) (61,1.937154589)
(61.25,1.949112303) (61.5,1.96140547) (61.75,1.973994928) (62,1.986839621) (62.25,1.999896734)
(62.5,2.013121853) (62.75,2.026469138) (63,2.039891507) (63.25,2.053340835) (63.5,2.066768155)
(63.75,2.080123879) (64,2.093358016) (64.25,2.106420403) (64.5,2.119260935) (64.75,2.131829797)
(65,2.144077695) (65.25,2.155956092) (65.5,2.167417433) (65.75,2.178415376) (66,2.188905013)
(66.25,2.198843089) (66.5,2.208188214) (66.75,2.216901071) (67,2.224944623) (67.25,2.232284311)
(67.5,2.238888252) (67.75,2.244727434) (68,2.249775907) (68.25,2.254010987) (68.5,2.257413448)
(68.75,2.259967727) (69,2.261662145) (69.25,2.262489128) (69.5,2.262445449) (69.75,2.261532495)
(70,2.259756553) (70.25,2.257129136) (70.5,2.253667336) (70.75,2.249394227) (71,2.24433932)
(71.25,2.238539076) (71.5,2.232037477) (71.75,2.22488668) (72,2.217147737) (72.25,2.20889138)
(72.5,2.20019889) (72.75,2.191162993) (73,2.181888773) (73.25,2.172494545) (73.5,2.163112586)
(73.75,2.153889612) (74,2.144986829) (74.25,2.136579324) (74.5,2.128854562) (74.75,2.122009667)
(75,2.116247288) (75.25,2.111769839) (75.5,2.108772212) (75.75,2.107433259) (76,2.10790675)
(76.25,2.110312776) (76.5,2.114730765) (76.75,2.121195163) (77,2.129694456) (77.25,2.140173668)
(77.5,2.152539824) (77.75,2.166669466) (78,2.182417086) (78.25,2.199623446) (78.5,2.218122989)
(78.75,2.237749922) (79,2.258342798) (79.25,2.279747683) (79.5,2.301820116) (79.75,2.324426106)
(80,2.34744242) (80.25,2.370756379) (80.5,2.394265347) (80.75,2.417876033) (81,2.441503713)
(81.25,2.465071433) (81.5,2.488509229) (81.75,2.511753397) (82,2.534745822) (82.25,2.557433378)
(82.5,2.579767383) (82.75,2.601703119) (83,2.623199421) (83.25,2.644218303) (83.5,2.66472464)
(83.75,2.684685894) (84,2.704071869) (84.25,2.722854504) (84.5,2.741007697) (84.75,2.758507138)
(85,2.775330186) (85.25,2.791455744) (85.5,2.806864164) (85.75,2.821537154) (86,2.83545771)
(86.25,2.84861004) (86.5,2.86097952) (86.75,2.872552639) (87,2.883316959) (87.25,2.893261074)
(87.5,2.902374587) (87.75,2.91064808) (88,2.918073088) (88.25,2.924642084) (88.5,2.930348461)
(88.75,2.93518652) (89,2.939151453) (89.25,2.942239342) (89.5,2.944447145) (89.75,2.945772692)
(90,2.946214686)
};
\addlegendentry{$q_1-1$}
\addplot+[no marks, thick, dashed] coordinates {
(0,-2.898704145) (0.25,-2.898171458) (0.5,-2.89657363) (0.75,-2.893911366) (1,-2.890185843)
(1.25,-2.885398704) (1.5,-2.879552061) (1.75,-2.872648499) (2,-2.864691062) (2.25,-2.855683267)
(2.5,-2.845629089) (2.75,-2.834532972) (3,-2.822399818) (3.25,-2.809234987) (3.5,-2.7950443)
(3.75,-2.779834029) (4,-2.763610902) (4.25,-2.746382095) (4.5,-2.728155233) (4.75,-2.708938381)
(5,-2.68874005) (5.25,-2.667569188) (5.5,-2.645435175) (5.75,-2.622347824) (6,-2.598317373)
(6.25,-2.573354485) (6.5,-2.547470241) (6.75,-2.520676136) (7,-2.492984077) (7.25,-2.464406375)
(7.5,-2.434955741) (7.75,-2.404645283) (8,-2.373488499) (8.25,-2.341499274) (8.5,-2.308691872)
(8.75,-2.27508093) (9,-2.240681458) (9.25,-2.205508824) (9.5,-2.169578757) (9.75,-2.132907337)
(10,-2.095510987) (10.25,-2.057406471) (10.5,-2.018610882) (10.75,-1.979141644) (11,-1.939016497)
(11.25,-1.898253493) (11.5,-1.856870993) (11.75,-1.814887656) (12,-1.772322432) (12.25,-1.72919456)
(12.5,-1.685523554) (12.75,-1.641329205) (13,-1.596631564) (13.25,-1.551450942) (13.5,-1.505807906)
(13.75,-1.459723263) (14,-1.41321806) (14.25,-1.366313577) (14.5,-1.319031321) (14.75,-1.271393016)
(15,-1.223420604) (15.25,-1.175136233) (15.5,-1.126562258) (15.75,-1.07772123) (16,-1.028635899)
(16.25,-0.9793292048) (16.5,-0.9298242759) (16.75,-0.8801444273) (17,-0.8303131602) (17.25,-0.7803541627)
(17.5,-0.7302913064) (17.75,-0.6801486523) (18,-0.6299504542) (18.25,-0.5797211628) (18.5,-0.5294854356)
(18.75,-0.4792681415) (19,-0.4290943804) (19.25,-0.3789894941) (19.5,-0.3289790851) (19.75,-0.2790890431)
(20,-0.2293455685) (20.25,-0.1797752083) (20.5,-0.130404897) (20.75,-0.08126200158) (21,-0.03237437739)
(21.25,0.01622956169) (21.5,0.06452077232) (21.75,0.1124694797) (22,0.1600450759) (22.25,0.2072159856)
(22.5,0.2539495148) (22.75,0.3002116688) (23,0.3459669341) (23.25,0.3911780238) (23.5,0.4358055713)
(23.75,0.4798077629) (24,0.5231399012) (24.25,0.565753881) (24.5,0.6075975574) (24.75,0.6486139858)
(25,0.6887405101) (25.25,0.7279076538) (25.5,0.7660377906) (25.75,0.8030435362) (26,0.8388258192)
(26.25,0.873271569) (26.5,0.9062509793) (26.75,0.9376142959) (27,0.9671881395) (27.25,0.9947714142)
(27.5,1.020130979) (27.75,1.042997479) (28,1.063061985) (28.25,1.079974534) (28.5,1.093346127)
(28.75,1.102756133) (29,1.107767138) (29.25,1.107948511) (29.5,1.102908115) (29.75,1.092328628)
(30,1.076001798) (30.25,1.053852511) (30.5,1.025946203) (30.75,0.9924778198) (31,0.953746025)
(31.25,0.9101199515) (31.5,0.8620060083) (31.75,0.8098201118) (32,0.7539676091) (32.25,0.6948308186)
(32.5,0.6327627392) (32.75,0.5680851501) (33,0.5010894568) (33.25,0.4320390605) (33.5,0.3611723975)
(33.75,0.2887061175) (34,0.2148381253) (34.25,0.1397503269) (34.5,0.06361102409) (34.75,-0.01342303393)
(35,-0.09120491825) (35.25,-0.169596081) (35.5,-0.2484651651) (35.75,-0.3276870206) (36,-0.4071418977)
(36.25,-0.4867147965) (36.5,-0.5662949246) (36.75,-0.6457752687) (37,-0.7250522475) (37.25,-0.8040254244)
(37.5,-0.8825972913) (37.75,-0.9606730879) (38,-1.038160671) (38.25,-1.114970401) (38.5,-1.191015072)
(38.75,-1.266209851) (39,-1.340472241) (39.25,-1.413722059) (39.5,-1.485881423) (39.75,-1.556874752)
(40,-1.626628777) (40.25,-1.695072546) (40.5,-1.762137454) (40.75,-1.827757257) (41,-1.891868102)
(41.25,-1.954408557) (41.5,-2.015319637) (41.75,-2.074544836) (42,-2.13203016) (42.25,-2.187724161)
(42.5,-2.241577956) (42.75,-2.29354527) (43,-2.343582459) (43.25,-2.391648537) (43.5,-2.437705199)
(43.75,-2.481716853) (44,-2.523650635) (44.25,-2.56347643) (44.5,-2.601166887) (44.75,-2.63669744)
(45,-2.670046314) (45.25,-2.701194539) (45.5,-2.730125952) (45.75,-2.756827202) (46,-2.781287754)
(46.25,-2.803499885) (46.5,-2.823458674) (46.75,-2.841162002) (47,-2.856610534) (47.25,-2.869807703)
(47.5,-2.880759698) (47.75,-2.889475434) (48,-2.895966538) (48.25,-2.900247306) (48.5,-2.902334684)
(48.75,-2.902248226) (49,-2.90001005) (49.25,-2.89564481) (49.5,-2.889179631) (49.75,-2.880644072)
(50,-2.870070069) (50.25,-2.857491868) (50.5,-2.842945979) (50.75,-2.826471094) (51,-2.808108026)
(51.25,-2.787899636) (51.5,-2.765890746) (51.75,-2.742128066) (52,-2.716660098) (52.25,-2.689537053)
(52.5,-2.66081075) (52.75,-2.630534515) (53,-2.598763079) (53.25,-2.565552464) (53.5,-2.530959874)
(53.75,-2.495043571) (54,-2.457862753) (54.25,-2.419477426) (54.5,-2.379948268) (54.75,-2.339336496)
(55,-2.297703718) (55.25,-2.255111796) (55.5,-2.211622684) (55.75,-2.167298283) (56,-2.122200282)
(56.25,-2.076389999) (56.5,-2.029928219) (56.75,-1.982875031) (57,-1.935289663) (57.25,-1.887230326)
(57.5,-1.83875404) (57.75,-1.789916486) (58,-1.740771842) (58.25,-1.691372635) (58.5,-1.641769592)
(58.75,-1.592011506) (59,-1.542145101) (59.25,-1.492214917) (59.5,-1.442263195) (59.75,-1.392329793)
(60,-1.342452089) (60.25,-1.292664927) (60.5,-1.243000564) (60.75,-1.193488629) (61,-1.144156121)
(61.25,-1.095027403) (61.5,-1.046124221) (61.75,-0.9974657547) (62,-0.949068662) (62.25,-0.9009471654)
(62.5,-0.8531131363) (62.75,-0.8055762072) (63,-0.7583438963) (63.25,-0.7114217357) (63.5,-0.6648134269)
(63.75,-0.6185209946) (64,-0.5725449563) (64.25,-0.5268844933) (64.5,-0.4815376309) (64.75,-0.4365014201)
(65,-0.391772125) (65.25,-0.3473454063) (65.5,-0.3032165085) (65.75,-0.2593804464) (66,-0.2158321871)
(66.25,-0.1725668284) (66.5,-0.1295797841) (66.75,-0.08686695253) (67,-0.0444248957) (67.25,-0.002251011322)
(67.5,0.039656299) (67.75,0.08129745481) (68,0.1226715227) (68.25,0.1637760391) (68.5,0.2046068283)
(68.75,0.2451578069) (69,0.2854207838) (69.25,0.3253852407) (69.5,0.365038091) (69.75,0.4043634266)
(70,0.4433422231) (70.25,0.4819520206) (70.5,0.5201665596) (70.75,0.5579553732) (71,0.595283322)
(71.25,0.6321100767) (71.5,0.6683895191) (71.75,0.7040690861) (72,0.7390890238) (72.25,0.7733815784)
(72.5,0.8068701157) (72.75,0.8394681896) (73,0.8710786) (73.25,0.9015925019) (73.5,0.9308886325)
(73.75,0.9588328004) (74,0.9852778095) (74.25,1.010064015) (74.5,1.033020802) (74.75,1.053969251)
(75,1.072726241) (75.25,1.089110152) (75.5,1.102948116) (75.75,1.114084482) (76,1.122389813)
(76.25,1.12776943) (76.5,1.130170343) (76.75,1.129585519) (77,1.126054809) (77.25,1.119662389)
(77.5,1.110531242) (77.75,1.098815591) (78,1.084692399) (78.25,1.068353003) (78.5,1.049995627)
(78.75,1.029819255) (79,1.008018982) (79.25,0.9847827852) (79.5,0.9602895062) (79.75,0.9347077929)
(80,0.9081957483) (80.25,0.8809010757) (80.5,0.8529615269) (80.75,0.8245055409) (81,0.7956529572)
(81.25,0.7665157608) (81.5,0.737198796) (81.75,0.7078004336) (82,0.6784131909) (82.25,0.6491242726)
(82.5,0.620016055) (82.75,0.5911665045) (83,0.5626495482) (83.25,0.5345353799) (83.5,0.5068907245)
(83.75,0.4797790695) (84,0.4532608505) (84.25,0.4273936094) (84.5,0.4022321269) (84.75,0.3778285333)
(85,0.3542323936) (85.25,0.3314907855) (85.5,0.3096483516) (85.75,0.2887473517) (86,0.2688276953)
(86.25,0.2499269687) (86.5,0.23208046) (86.75,0.2153211693) (87,0.1996798215) (87.25,0.1851848725)
(87.5,0.1718625123) (87.75,0.159736664) (88,0.1488289871) (88.25,0.1391588722) (88.5,0.1307434417)
(88.75,0.1235975407) (89,0.1177337407) (89.25,0.1131623304) (89.5,0.1098913157) (89.75,0.1079264183)
(90,0.1072710705)
};
\addlegendentry{$q_2-1$}
\end{axis}
\end{tikzpicture}
\caption{Signed phase-slowness errors of the two Maxwell IOFD polarizations
at $\ppw=6$, with fixed azimuth $\beta=\pi/7$ and
$\dvec=(\sin\alpha\cos\beta,\sin\alpha\sin\beta,\cos\alpha)$.
The two branches coincide on the coordinate axis and separate away from
it.  Values are computed from the radial roots.}
\label{fig:angular6}
\end{minipage}
\end{figure}

\subsection{Source and field filter}\label{sec:amplitude}

The coefficient fit places the numerical wave surfaces; the amplitudes of
the waves radiated by a point source are set by the residues of the
discrete Green symbol at these surfaces (Proposition~\ref{prop:residue}).
We move these residues
toward those of the continuous symbol with a compact filter $Q_h$ \cite[Sections~2.4 and~5.2]{Stolk2016}, applied
identically to each component with weights $Q_m$ on the shells
$\mathcal N_m=\{\vct\delta\in\{-1,0,1\}^3:|\vct\delta|_1=m\}$,
\begin{equation}
 Q_0=1-3\beta_1-3\beta_2-\beta_3,\qquad
 Q_1=\frac{\beta_1}{2},\quad Q_2=\frac{\beta_2}{4},\quad
 Q_3=\frac{\beta_3}{8}.
 \label{eq:Qweights}
\end{equation}
For a discrete source load $\vct f_h$, which assigns $h^{-3}$ times the
source strength to its selected electric component, we solve and reconstruct
\begin{equation}
 A_h\vct v_h=Q_h\vct f_h=: \vct b_h,\qquad
 \E_h=Q_h\vct v_h.
 \label{eq:fullscheme}
\end{equation}
In a homogeneous medium this solution operator has the symbol
$h^2q_Q^2\mathscr P^{-1}$, so the filter changes the modal residues but
not the dispersion surfaces. For the resolution of 6 points per
wavelength used in section~\ref{sec:field-validation},
\begin{equation}
 (\beta_1,\beta_2,\beta_3)
 =(-0.2002664857,\ 0.1236538386,\ 0.0484221099).
 \label{eq:Qcoeffs}
\end{equation}

\subsection{Variable materials and perfectly matched layers}\label{sec:pml}

To apply the filtered Maxwell scheme \eqref{eq:fullscheme} in a variable
medium, we select its face and edge coefficients from the local material
values. Freezing $\epsilon$ and $\mu$ at a point gives the normalized
equation \eqref{eq:pde} with wavenumber $k=k_0\sqrt{\epsilon\mu}$.
For each face or edge row, we therefore evaluate the interpolation
\eqref{eq:linearinterpolation} at the local dimensionless wavenumber
\begin{equation}
 \eta(\x)=\frac{k_0h\sqrt{\epsilon(\x)\mu(\x)}}{2\pi}.
 \label{eq:localeta}
\end{equation}
Here $\x$ is the row center on the component's own staggered lattice.
This selection applies the frozen-material dispersion design at the
local resolution, while the sampled material values enter the
constitutive factors of the Maxwell operator.

To impose the outgoing condition for this variable-material problem on
a finite domain, we surround $\Omega_{\ph}$ with a Cartesian PML of
thickness $L_{\pml}$ and denote the resulting box by $\Omega$. Let $d_j^{\pml}(\x)$ be the distance from $x_j$ to the
physical interval in direction $j$, so that $d_j^{\pml}=0$ inside
$\Omega_{\ph}$. We use the stretching factors
\begin{equation}
 \zeta_j=1+\ii\rho_j,\qquad
 \rho_j=\frac{2c_{\pml}}{k_0L_{\pml}}
                    \frac{d_j^{\pml}}{L_{\pml}},\qquad
 T=\diag\!\left(\frac{\zeta_2\zeta_3}{\zeta_1},
                \frac{\zeta_1\zeta_3}{\zeta_2},
                \frac{\zeta_1\zeta_2}{\zeta_3}\right),
 \label{eq:stretch}
\end{equation}
whose absorption profile $\rho_j$ grows linearly with the depth into the
layer, so that the complex coordinate stretching is quadratic. The
stretching gives
\begin{equation}
 \curl(\mu^{-1}T^{-1}\curl\E)-k_0^2\epsilon T\E=\f
       \quad\text{in }\Omega.
 \label{eq:pmlpde}
\end{equation}
This constitutive PML formulation \cite{Berenger1994,ChewWeedon1994}
requires material values throughout the layer. We obtain them by normal
extension: each coordinate is projected onto its physical interval before
$\epsilon$ and $\mu$ are evaluated. This retains the tangential material
variation and continues the boundary values along each outward normal.
On the exterior boundary of the PML, we impose a first-order
Silver--M\"uller absorbing condition in stretched coordinates.

Inside the PML, the real local wavenumber \eqref{eq:localeta} still
selects the IOFD coefficients, and the stretching enters only through
the constitutive factors. We represent these factors by diagonal face and
edge matrices with staggered samples
\begin{equation}
 (D_{\mu,h})_r=\left(\mu^{-1}\frac{\zeta_r}{\zeta_s\zeta_t}\right)^{1/2},
 \qquad
 (D_{\epsilon,h})_r
       =\left(\epsilon\frac{\zeta_s\zeta_t}{\zeta_r}\right)^{1/2},
 \label{eq:sqrtfactors}
\end{equation}
where the complex square roots are continued from their positive values
in the physical medium. The variable-coefficient operator is
\begin{equation}
 A_h=C_h^*\mathsf F_hC_h-k_0^2\mathsf M_h,\qquad
 \mathsf F_h=D_{\mu,h}F_h^\Omega D_{\mu,h},\quad
 \mathsf M_h=D_{\epsilon,h}M_h^\Omega D_{\epsilon,h}.
 \label{eq:pmloperator}
\end{equation}
Here $F_h^\Omega$ and $M_h^\Omega$ are the locally selected filters on
the truncated lattices. Rows whose filter support lies in the lattice use
\eqref{eq:face}--\eqref{eq:edge}; at the outer boundary, filter rows with
incomplete support are replaced by the identity, and the tangential
boundary rows discretize the absorbing condition. The square-root factors
enter on both sides without complex conjugation, so that
$\mathsf F_h$ and $\mathsf M_h$ reproduce the complex constitutive
tensors of \eqref{eq:pmlpde}.

With \eqref{eq:pmloperator}, the filtered system \eqref{eq:fullscheme}
now represents the outgoing problem on a finite domain. Its compact
structure also permits the same construction on the coarser grids used
by the solver below.

\section{Three-Grid Preconditioner}\label{sec:threegrid}

In this section, we construct
a three-grid preconditioner for the IOFD system that assigns
two distinct roles to its coarse levels.
An unshifted intermediate-grid IOFD equation supplies a global correction
whose propagation is governed by the optimized coarse-grid dispersion
relation. A complex-shifted auxiliary cycle introduces damping to make
the iterative solution of that equation inexpensive. Since every level uses the
discretization of Section~\ref{sec:discretization}, all operators remain
compact. We first define the level operators and transfers
(Section~\ref{sec:levels}), then the nested iteration
(Section~\ref{sec:cycle}), and finally its cost and distributed
implementation (Section~\ref{sec:implementation}).

\subsection{Level operators and edge transfers}\label{sec:levels}

To carry the local dispersion design into the hierarchy, we discretize
the operator \eqref{eq:pmloperator} at each level's own grid spacing and
material samples. Let $X_h$, $X_{2h}$, and $X_{4h}$ be the
edge spaces on three nested grids covering the same box $\Omega$. The
spacing doubles and the number of PML cells halves from one level to the
next, so the PML has the same physical width on every level. On a level
$\ell\in\{h,2h,4h\}$, write
\begin{equation}
 K_\ell=C_\ell^*\mathsf F_\ell C_\ell,\qquad
 A_\ell=K_\ell-k_0^2\mathsf M_\ell,
 \label{eq:physicalpair}
\end{equation}
where the filters in $\mathsf F_\ell$ and $\mathsf M_\ell$ are selected
by the local dimensionless wavenumber \eqref{eq:localeta} with $h$
replaced by $\ell$. The physical
intermediate-grid equation uses $A_{2h}$. For heterogeneous media, the
$2h$ operators use, at the knots $0.175\le\eta_j\le0.325$, coefficients
obtained by continuation from the primary coefficients at $\eta=1/3$
toward smaller $\eta$ under the dispersion and filter-positivity
constraints; this extends the small-$a$ coefficient family down to
$\eta=0.175$, and the same linear interpolation applies. For these heterogeneous $2h$
operators, on a boundary face of $\Omega$ normal to the $x_r$-direction,
the mass row of the normal component $E_r$ also retains its transverse
filter $I+b_t(\Delta_s+\Delta_t)+b_{tt}\Delta_s\Delta_t$ wherever the
transverse support is complete.

The correction on a coarser grid is transferred using the staggered
geometry of each electric component. On field-component coordinates,
the one-dimensional nodal and edge prolongations are
\begin{align}
 (P_nw)_{2j}=w_j,&\qquad (P_nw)_{2j+1}=\tfrac12(w_j+w_{j+1}),\\
 (P_ew)_{2j}=w_j,&\qquad (P_ew)_{2j+1}=w_j.
\end{align}
Tensorization gives
\begin{equation}
 P_\ell=\diag(P_e\otimes P_n\otimes P_n,
               P_n\otimes P_e\otimes P_n,
               P_n\otimes P_n\otimes P_e),\qquad
 R_\ell=\tfrac18P_\ell^*,\quad \ell=h,2h.
 \label{eq:transfers}
\end{equation}
The factor $1/8$ corresponds to the pointwise finite-difference scaling
of the operators. These geometric transfers copy an edge component in
its own direction and interpolate linearly in the transverse directions.

\subsection{The auxiliary cycle and physical correction}\label{sec:cycle}

To compute the correction transferred by $P_h$, we solve the physical
$A_{2h}$ equation iteratively. We precondition this solve with an
auxiliary two-grid cycle whose operators on $2h$ and $4h$ include a
complex shift. For a shift parameter $\sigma$, define
\begin{equation}
 B_\ell=K_\ell-(1+\ii\sigma)k_0^2\mathsf M_\ell,
 \qquad \ell=2h,4h,\quad \sigma>0.
 \label{eq:auxiliary}
\end{equation}
The shift multiplies only the mass term: the filters are those of the
unshifted operator on the same level, selected by the real local
wavenumber, and the impedance of the absorbing boundary condition is
scaled consistently by $\sqrt{1+\ii\sigma}$. The damping makes the
auxiliary equations amenable to Jacobi-preconditioned Krylov smoothing
and coarse correction, and the resulting cycle serves as an approximate
inverse of the physical intermediate-grid operator
\cite{ErlanggaEtAl2006,FuEtAlHelmholtz,FuMaxwell}.

All inner solves perform a prescribed amount of work. Let
$\calS_H^{(m)}(g;x_0)$ denote $m$ GMRES steps for an operator $H$,
applied to $Hx=g$ with initial guess $x_0$. Each step is right-preconditioned
by 2 damped Jacobi sweeps. For an input $g$, these sweeps start from
$y^{(0)}=0$ and use
\begin{equation}
 y^{(j+1)}=y^{(j)}+\omega_J\diag(H)^{-1}(g-Hy^{(j)}),
 \qquad j=0,1.
 \label{eq:Jacobi}
\end{equation}
The auxiliary cycle combines smoothing on $2h$, a correction on $4h$,
and further smoothing on $2h$. Its action on a right-hand side $g$ is
\begin{equation}
 \begin{aligned}
 y^{\rm pre}&=\calS_{B_{2h}}^{(m_s)}(g;0),\\
 g_{4h}&=R_{2h}(g-B_{2h}y^{\rm pre}),\\
 e_{4h}&=\calS_{B_{4h}}^{(m_b)}(g_{4h};0),\\
 \calV_{2h}^{\sigma}(g)
   &=\calS_{B_{2h}}^{(m_s)}(g;y^{\rm pre}+P_{2h}e_{4h}).
 \end{aligned}
 \label{eq:Vcycle}
\end{equation}
The physical intermediate-grid equation is approximately solved by
$m_c$ FGMRES steps, with $\calV_{2h}^{\sigma}$ as a right preconditioner
and zero initial guess. Let $\calK_{2h}(r_{2h})$ denote the resulting
approximation to $A_{2h}^{-1}r_{2h}$. With
$\calS_h=\calS_{A_h}^{(m_h)}$, Algorithm~\ref{alg:3g} combines this
correction with fine-grid pre- and post-smoothing.

\begin{algorithm}[H]
\caption{Three-grid preconditioner $z_h=\calB_h(r_h)$}
\label{alg:3g}
\begin{algorithmic}[1]
\State $z_h^{\rm pre}\gets\calS_h(r_h;0)$.
\State $r_{2h}\gets R_h(r_h-A_hz_h^{\rm pre})$.
\State $e_{2h}\gets\calK_{2h}(r_{2h})$ using $m_c$ FGMRES steps on $A_{2h}$,
       preconditioned by $\calV_{2h}^{\sigma}$.
\State $\widetilde z_h\gets z_h^{\rm pre}+P_he_{2h}$.
\State $z_h\gets\calS_h(r_h;\widetilde z_h)$.
\end{algorithmic}
\end{algorithm}

Because the inner Krylov iterations depend nonlinearly on their input
residuals, $\calB_h$ is a variable preconditioner. We therefore apply
restarted FGMRES \cite{Saad1993} to the fine-grid system
\eqref{eq:fullscheme} and stop when the directly evaluated relative
residual satisfies
\begin{equation}
 \frac{\norm{A_h\vct v_h-\vct b_h}_2}{\norm{\vct b_h}_2}
       \le\varepsilon_{\rm out}.
 \label{eq:stop}
\end{equation}
The tolerance $\varepsilon_{\rm out}$ is chosen in
Section~\ref{sec:field-validation} by comparison with an analytic
wavefield.

\subsection{Cost per application and distributed implementation}
\label{sec:implementation}

We quantify the cost of Algorithm~\ref{alg:3g} by counting its
prescribed Krylov steps on each level. Each step uses compact stencil
evaluations, diagonal Jacobi scaling, and vector operations; geometric
transfers connect the levels. For fixed Krylov dimensions, the work and
storage are proportional to the number of fine-grid unknowns, and every
level is applied matrix-free. Specifically, one application performs $2m_h$ GMRES steps with $A_h$,
$m_c$ FGMRES steps with $A_{2h}$, and $m_c$ auxiliary cycles, each with
$2m_s$ GMRES steps with $B_{2h}$ and $m_b$ GMRES steps with $B_{4h}$.
For the values $(m_h,m_c,m_s,m_b)=(6,6,4,10)$ used in
Section~\ref{sec:gpu-results}, this gives 12 steps on the fine grid and
$6+48+60=114$ steps on the coarser grids. Although the $2h$ and $4h$
grids carry about $1/8$ and $1/64$ as many unknowns as the fine grid,
they account for most Krylov steps. Interface exchanges and global
inner products on these levels therefore remain important to the cost
of a distributed application.

The distributed GPU implementation partitions $\Omega$ in the $y$ and
$z$ directions, with each rank owning an $x$-spanning pencil, and the CPU
implementation partitions it in all three directions. The
partition is nested across the three grids, so geometric transfers
involve only neighboring pencils. Stencil evaluations exchange interface
data while independent interior rows are evaluated, and Krylov
orthogonalization uses global reductions. On the GPU platform of
Section~\ref{sec:gpu-results}, interface data are exchanged within a node
through the vendor collective library and between nodes through MPI.
The smoothing GMRES steps on all three grids are carried out in $s$-step
form, with a monomial Krylov basis and a single Gram reduction per
smoothing call. On the $h$ and $2h$ grids, the basis vectors after the
first are stored in bfloat16, while all arithmetic remains in single
precision; this halves their storage and memory traffic.

\section{Wavefield and Solver Validation}\label{sec:field-validation}

In this section, we compare the wavefields computed with the IOFD
discretization and the three-grid solver with the outgoing Maxwell Green
tensor, which provides an analytic reference for both phase and amplitude
over long propagation distances and includes the effects of the source
and field filter and the PML. We first measure how the field error grows
with the electrical size (Section~\ref{sec:green}) and then use the same
reference to select the stopping tolerance and the storage precision of
the smoother bases (Section~\ref{sec:tolerance}).

The source is a $z$-oriented point load on the Yee edge nearest to
$(0,0.5,0)$, filtered as in \eqref{eq:fullscheme}; its grid-edge location
$\x_s$ also locates the Green reference. We use 6 fine-grid points per wavelength,
12 PML cells on each side with strength $c_{\pml}=3$, and the first-order
absorbing condition on the exterior of the PML. A source-centered ball
$B_{\lambda}$ of radius one wavelength, $\lambda=2\pi/k$, is excluded
from the field comparisons \cite{TournierEtAl2022}.

For the normalized homogeneous equation, the outgoing Green tensor is
\begin{equation}
 \mathsf G_k(\x)=\left(I+k^{-2}\nabla\nabla^T\right)
                       \frac{\exp(\ii k|\x|)}{4\pi|\x|}.
 \label{eq:greentensor}
\end{equation}
A unit $z$-oriented point load gives the field
$\mathsf G_k(\x-\x_s)\vct e_3$.

We compare the computed field $E_{z,h}$ with
$E_z^{\rm ref}=G_{zz}(\x-\x_s)$ directly at the native $z$-edge locations.
With
$r_j=|\bm{x}^{(j)}-\bm{x}_s|$ and
$\mathcal I$ denoting the $z$-edges in $\Omega_{\rm ph}\setminus B_{\lambda}$,
the field error is
\begin{equation}\label{eq:iofd-green-err}
{\rm Err}=
\frac{\sum_{j\in\mathcal I}r_j|\operatorname{Re}(E_{z,h}^{(j)}-E_z^{{\rm ref},(j)})|}
     {\sum_{j\in\mathcal I}r_j|\operatorname{Re}E_z^{{\rm ref},(j)}|}
+
\frac{\sum_{j\in\mathcal I}r_j|\operatorname{Im}(E_{z,h}^{(j)}-E_z^{{\rm ref},(j)})|}
     {\sum_{j\in\mathcal I}r_j|\operatorname{Im}E_z^{{\rm ref},(j)}|}.
\end{equation}
We also report the relative vector $L^2$ error using all three field
components on their respective edge lattices in the same comparison region.

\subsection{Comparison with the Green field}\label{sec:green}

To examine how propagation errors accumulate at fixed local resolution,
we increase the mesh from $480^3$ to $5760^3$ at 6 points per wavelength,
so that the physical domain spans 76 to 956 wavelengths in each
coordinate direction. The meshes up to $2880^3$ are solved on the GPUs
of Section~\ref{sec:gpu-results}, and the $3840^3$ and $5760^3$ meshes with the CPU
implementation of Section~\ref{sec:cpu}. The solves use the
relative residual tolerance $\varepsilon_{\rm out}=10^{-3}$, justified in
Section~\ref{sec:tolerance}.
Table~\ref{tab:iofd-green} summarizes the errors.

\begin{table}[htbp]
\centering\small
\caption{Homogeneous Green-field errors. Mesh dimensions include the PML.}
\label{tab:iofd-green}
\begin{tabular}{@{}ccccc@{}}
\toprule
mesh & wavelengths & platform & Err & vector $L^2$ \\
\midrule
$480^3$ & 76 & GPU & $2.8326\%$ & $1.3560\%$ \\
$960^3$ & 156 & GPU & $2.9496\%$ & $1.3840\%$ \\
$1920^3$ & 316 & GPU & $4.4813\%$ & $1.8921\%$ \\
$2880^3$ & 476 & GPU & $5.7995\%$ & $2.4474\%$ \\
$3840^3$ & 636 & CPU & $7.1925\%$ & $3.0586\%$ \\
$5760^3$ & 956 & CPU & $10.1930\%$ & $4.3765\%$ \\
\bottomrule
\end{tabular}
\end{table}

Over this more than twelvefold increase in electrical size at fixed
6-point resolution, the vector $L^2$ error rises only from $1.36\%$ to
$4.38\%$ and Err from $2.83\%$ to $10.19\%$. The slow growth reflects the
dispersion design: by \eqref{eq:accumulation}, the sampled wavenumber
errors of Table~\ref{tab:dispersion} accumulate over 956 wavelengths to a
phase discrepancy of about
$2\pi\cdot956\cdot3.08\times10^{-5}\approx0.19$ radians for IOFD,
compared with 315 radians for standard Yee.

\begin{figure}[htbp]
\centering
\includegraphics[width=\textwidth]{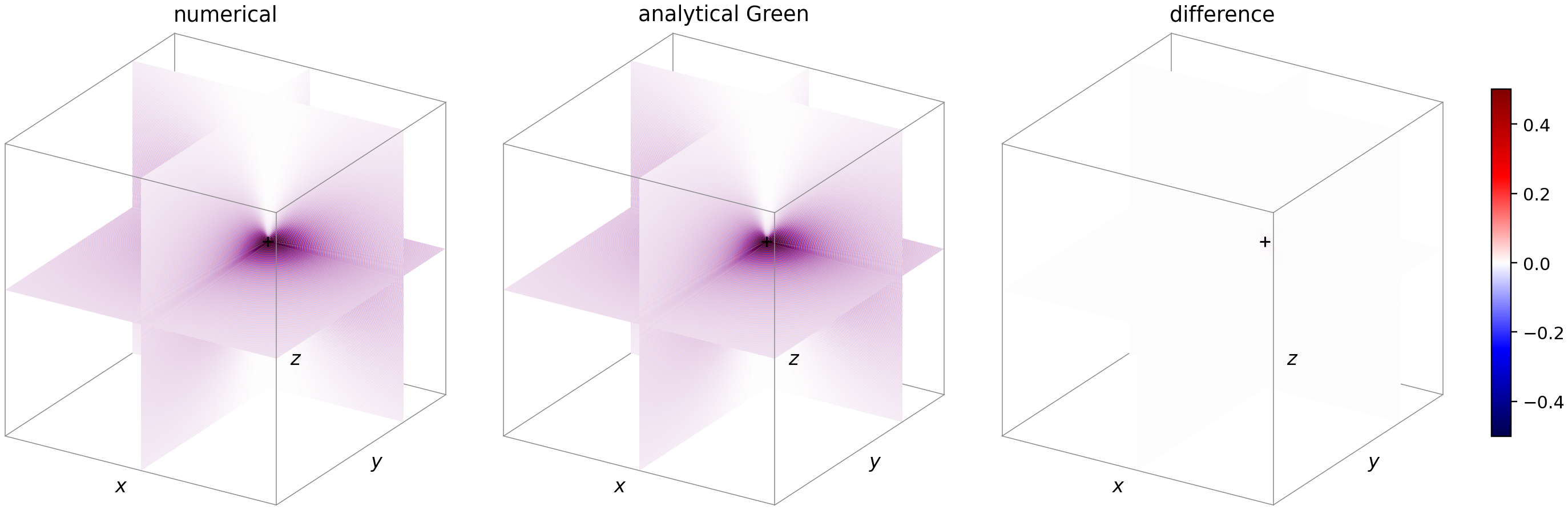}
\caption{Computed field, outgoing Green field, and their difference on the
$x$-, $y$-, and $z$-normal planes through the source for the homogeneous
$960^3$ problem, shown in the physical domain $[-1,1]^3$.  The plotted
quantity is $\operatorname{Re}(E_z)$; the two fields are sampled at identical
native $z$-edge locations.  The cross marks the discrete source, and the
difference is set to zero within the exclusion radius in
Eq.~\eqref{eq:iofd-green-err}.  All three panels share one color scale.}
\label{fig:iofd-green-slices}
\end{figure}

\begin{figure}[htbp]
\centering
\includegraphics[width=0.8\textwidth]{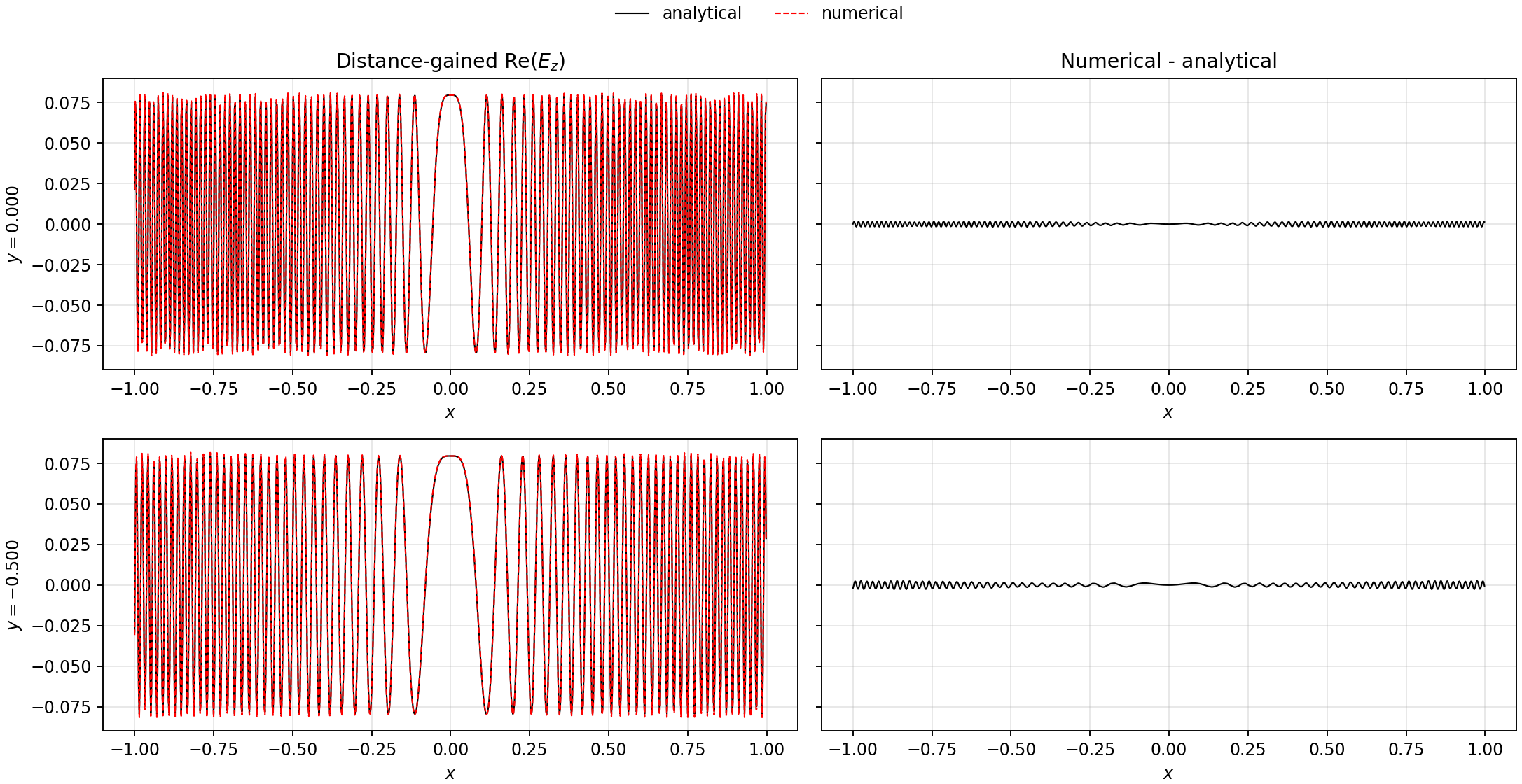}
\caption{One-dimensional receiver-line comparisons on the $960^3$ homogeneous
problem.  The left panels compare $r\operatorname{Re}(E_z)$ with the outgoing
Green reference; the right panels show the numerical minus analytical
difference on the same lines.}
\label{fig:iofd-green-lines}
\end{figure}

The same accuracy holds throughout the domain. On the $960^3$ mesh solved
on 8 GPUs, Figure~\ref{fig:iofd-green-slices} displays the computed and
analytic $E_z$ on the three planes through the source, together with
their difference on the same color scale. Because the field decays with
the distance from the source, Figure~\ref{fig:iofd-green-lines} compares
distance-weighted receiver profiles along two lines away from the source,
where the computed field follows the phase and amplitude of the Green
field.

\subsection{Stopping residual and basis storage}\label{sec:tolerance}

All computations use complex single-precision arithmetic, and the
smoother bases on the $h$ and $2h$ grids are stored in bfloat16
(Section~\ref{sec:implementation}). We choose the stopping tolerance by
relating the algebraic residual to the error against the analytic field,
and we check the reduced basis storage at the same time. For the
homogeneous $960^3$ problem on 8 GPUs, Table~\ref{tab:iofd-tolerance}
reports the field error and the number of preconditioner applications as
the residual target decreases from $10^{-2}$ to $10^{-5}$, with the
smoother bases stored in bfloat16 and in single precision.

\begin{table}[htbp]
\centering\small
\caption{Effect of the stopping residual on field accuracy for the homogeneous $960^3$ problem, with the smoother bases stored in bfloat16 / single precision.}
\label{tab:iofd-tolerance}
\begin{tabular}{@{}cccc@{}}
\toprule
outer target & PC calls & $\|A_hv_h-b_h\|_2/\|b_h\|_2$ & vector $L^2$ error \\
\midrule
$10^{-2}$ & 14 / 14 & $9.52\times10^{-3}$ / $8.94\times10^{-3}$ & 21.28\% / 18.86\% \\
$10^{-3}$ & 23 / 22 & $6.17\times10^{-4}$ / $9.50\times10^{-4}$ & 1.384\% / 1.445\% \\
$10^{-4}$ & 28 / 27 & $9.25\times10^{-5}$ / $9.78\times10^{-5}$ & 1.350\% / 1.356\% \\
$10^{-5}$ & 36 / 35 & $8.48\times10^{-6}$ / $9.08\times10^{-6}$ & 1.344\% / 1.344\% \\
\bottomrule
\end{tabular}
\end{table}

\FloatBarrier
Tightening the target from $10^{-2}$ to $10^{-3}$ reduces the field
error from about $20\%$ to $1.4\%$. Further residual reduction has little
effect on field accuracy: a target of $10^{-5}$ lowers the error only to
$1.34\%$, the discretization error of this mesh, but requires 36 instead
of 23 preconditioner applications. The bfloat16 basis storage changes
the number of applications by at most one, and at the tightest target
both variants reach the same field error, so the reduced storage does
not limit the attainable accuracy. We therefore use
$\varepsilon_{\rm out}=10^{-3}$ in all subsequent experiments.

\FloatBarrier

\section{Numerical Results}\label{sec:gpu-results}

In this section, we evaluate the efficiency and parallel scalability of
the three-grid solver on GPUs and CPUs, using the precision and stopping
criterion established in Section~\ref{sec:field-validation}. Weak-scaling
tests measure how the cost grows with the domain size in wavelengths at
a fixed number of points per wavelength, while strong-scaling tests
measure how additional resources reduce the solution time for a fixed
problem. After specifying the media and common settings
(Section~\ref{sec:settings}), we report weak and strong scaling on GPUs
(Sections~\ref{sec:weak} and~\ref{sec:strong}) and then extend the
computation to 574 billion unknowns on CPUs (Section~\ref{sec:cpu}).

\subsection{Media and common settings}\label{sec:settings}

We take four media from the Maxwell benchmarks of \cite{FuMaxwell},
which test the three-grid method under uniform
propagation, smooth focusing, periodic variation, and spatial disorder.
The physical domain is $(-1,1)^3$. All heterogeneous media use
$\epsilon(\x)=\mu(\x)=m(\x)=1/c(\x)$, with material values extended along the outward normals into the PML as
specified in Section~\ref{sec:pml}.
The material fields are defined as follows:
\begin{enumerate}
\item \emph{Homogeneous medium}: $m=1$.
\item \emph{Converging lens}:
\begin{equation}
 m(\x)=\frac{3/4}{1-\tfrac12\exp(-8|\x|^2)}.
 \label{eq:lensmodel}
\end{equation}
\item \emph{Periodic medium}:
\begin{equation}
 m(\x)=\left[1+\frac14\left\{
 \cos\!\left(20\frac{x+y}{\sqrt2}\right)
 +\cos\!\left(20\frac{x-y}{\sqrt2}\right)\right\}\right]^{1/2}.
 \label{eq:periodicmodel}
\end{equation}
\item \emph{Random medium}: a fixed 8-term Fourier field
\begin{equation}
 m(\x)=\frac{4/3}{1-\tfrac12\delta(\x)},\qquad
 \delta(\x)=\frac12+\frac{0.45}{3.32}R(\x),
 \label{eq:randommodel}
\end{equation}
where $R$ is specified in Appendix~\ref{app:random}.
\end{enumerate}

Every mesh size $n^3$ counts fine-grid intervals including 12 PML cells
on each side. With $h=2/(n-24)$, the frequency is set by
$k_0hm_*=2\pi/6$, where the reference maximum slowness $m_*$ is
$1$, $1.5$, $\sqrt{1.5}$, and $2.32675$ for the homogeneous,
lens, periodic, and random media, respectively. Increasing $n$
therefore increases the electrical size of the fixed physical domain
while retaining the same minimum number of points per wavelength. All sources use the $z$-oriented grid-edge load described in
Section~\ref{sec:field-validation}, with the source-and-field filter of
\eqref{eq:fullscheme}.

The PML uses the linear absorption profile \eqref{eq:stretch} and the
first-order absorbing condition on its exterior. Its strength is tied to
the same wavelength normalization,
\begin{equation}
 c_{\pml}=\frac{3n_{\pml}}{2\ppw m_*},\qquad
 n_{\pml}=12,\quad\ppw=6.
 \label{eq:pmlstrength}
\end{equation}
For the homogeneous medium, this gives the value $c_{\pml}=3$ of
Section~\ref{sec:field-validation}. The auxiliary shift is $\sigma=0.2$
for the homogeneous medium and $\sigma=0.6$ for the heterogeneous media.
The Krylov work counts of Section~\ref{sec:cycle} are
\begin{equation}
 (m_h,m_c,m_s,m_b)=(6,6,4,10).
 \label{eq:workcounts}
\end{equation}
The Jacobi weights on $h$, $2h$, and $4h$ are $0.5$, $0.55$, and
$0.4$, respectively. The outer FGMRES restart length is 4, and all
reported solves satisfy \eqref{eq:stop} with
$\varepsilon_{\rm out}=10^{-3}$ in complex single precision.

One preconditioner application (PC) denotes one complete call of
Algorithm~\ref{alg:3g}. We report setup and complete solution times separately,
the cost per application $\tau=T_{\rm solve}/n_{\rm PC}$, and the
maximum device-memory usage per GPU.
All GPU runs use a cluster with eight Hygon BW1000 GPUs and four 400~Gb/s
network adapters per node. The BW1000 is a general-purpose GPU with
64~GB of HBM, a memory bandwidth of 1.8~TB/s, a peak single-precision
(FP32) throughput of 60~TFLOPS, a PCIe~5.0~$\times$16 host interface, and
a 448~GB/s bidirectional GPU interconnect. Each MPI process controls one GPU. Within one node,
the nested $y$--$z$ partition of Section~\ref{sec:implementation} uses
$z$-slabs only; multi-node runs use eight processes per node, and each
node owns one $z$-slab, divided into eight $y$-blocks.

\subsection{Weak scaling}\label{sec:weak}

We first increase the electrical size of the problem together with the
GPU allocation, keeping the number of unknowns per GPU nearly constant.
The sequence uses $480^3$ on 1 GPU, $960^3$ on 8 GPUs, $1920^3$ on
64 GPUs, and $2880^3$ on 216 GPUs. Under the frequency normalization
above, the homogeneous domain grows from 76 to 156, 316, and then 476
wavelengths in each coordinate direction. At the largest scale, the
solver must therefore resolve three-dimensional propagation across
nearly 480 wavelengths using approximately 71.7 billion complex edge
unknowns, within an essentially fixed storage budget per GPU.

Let $D(n)=3n(n+1)^2$ be the number of global Yee edge unknowns for a total
$n^3$ mesh, and let $\tau_{n,p}$ be the mean full-solve time per PC call on
$p$ GPUs. We correct for the small
change in the number of unknowns per GPU \cite{FuEtAlHelmholtz} and define
\begin{equation}
 \eta_{\rm w}(n,p)
 =\frac{\tau_{480,1}}{\tau_{n,p}}
  \frac{D(n)/p}{D(480)},
 \qquad (n,p)\in\{(480,1),(960,8),(1920,64),(2880,216)\}.
 \label{eq:iofd-weak-efficiency}
\end{equation}
Thus $\eta_{\rm w}$ measures the parallel throughput of one
application, whereas the PC count reflects the convergence of the solver
as the propagation distance grows. Table~\ref{tab:iofd-weak} reports both
and their combined effect on the solution time.

\begin{table}[htbp]
\centering
\caption{Weak scaling on GPUs, relative to the $480^3$ 1-GPU case for each medium. Memory is the maximum measured device usage per GPU.}
\label{tab:iofd-weak}
\footnotesize
\begin{adjustbox}{max width=\textwidth}
\begin{tabular}{ccccccccc}
\toprule
medium & mesh & PC calls & GPUs & $T_{\rm setup}$ & $T_{\rm solve}$ & solve/PC & $\eta_{\rm w}$ & peak/GPU \\
\midrule
homogeneous & $480^3$ & 14 & 1 & 0.240~s & 15.895~s & 1.1353~s & 100.0\% & 50.35~GiB \\
homogeneous & $960^3$ & 23 & 8 & 0.327~s & 24.706~s & 1.0742~s & 105.5\% & 51.91~GiB \\
homogeneous & $1920^3$ & 42 & 64 & 4.922~s & 49.375~s & 1.1756~s & 96.3\% & 46.83~GiB \\
homogeneous & $2880^3$ & 61 & 216 & 4.784~s & 77.480~s & 1.2702~s & 89.1\% & 48.57~GiB \\
\addlinespace
lens & $480^3$ & 22 & 1 & 2.335~s & 52.702~s & 2.3955~s & 100.0\% & 48.58~GiB \\
lens & $960^3$ & 45 & 8 & 2.673~s & 98.500~s & 2.1889~s & 109.2\% & 50.15~GiB \\
lens & $1920^3$ & 106 & 64 & 4.708~s & 239.432~s & 2.2588~s & 105.7\% & 50.25~GiB \\
lens & $2880^3$ & 174 & 216 & 4.704~s & 459.206~s & 2.6391~s & 90.5\% & 50.98~GiB \\
\addlinespace
periodic & $480^3$ & 26 & 1 & 2.330~s & 62.929~s & 2.4203~s & 100.0\% & 48.58~GiB \\
periodic & $960^3$ & 50 & 8 & 3.181~s & 110.465~s & 2.2093~s & 109.3\% & 50.15~GiB \\
periodic & $1920^3$ & 113 & 64 & 4.455~s & 257.479~s & 2.2786~s & 105.9\% & 50.25~GiB \\
periodic & $2880^3$ & 190 & 216 & 4.996~s & 503.773~s & 2.6514~s & 91.0\% & 50.98~GiB \\
\addlinespace
random & $480^3$ & 22 & 1 & 2.372~s & 53.242~s & 2.4201~s & 100.0\% & 48.58~GiB \\
random & $960^3$ & 41 & 8 & 3.747~s & 90.591~s & 2.2095~s & 109.3\% & 50.15~GiB \\
random & $1920^3$ & 84 & 64 & 4.382~s & 191.471~s & 2.2794~s & 105.8\% & 50.25~GiB \\
random & $2880^3$ & 136 & 216 & 4.671~s & 360.639~s & 2.6518~s & 90.9\% & 50.98~GiB \\
\bottomrule
\end{tabular}%
\end{adjustbox}
\end{table}

\FloatBarrier
The homogeneous problem spanning 476 wavelengths per direction converges
in 61 PC calls and 77.5 seconds. Relative to the 1-GPU case, the propagation
distance increases by a factor of 6.3, whereas the PC count grows by a
factor of about 4.4. The same three-grid construction solves every
heterogeneous $2880^3$ problem within 9 minutes: the random, lens, and
periodic media require 360.6--503.8 seconds.

The cost per application changes little from 1 to 64 GPUs, with
weak-scaling efficiencies of 96.3--109.3\% across the four media. At
216 GPUs, one application to a system with approximately $71.7$ billion
complex unknowns takes 1.27--2.65 seconds, retaining 89.1--91.0\%
weak-scaling efficiency. The matrix-free construction also keeps setup
and storage small: setup takes at most 5.0 seconds throughout the
sequence, and device-memory usage remains below 52 GiB per GPU as the
global mesh grows by a factor of 216.

\subsection{Strong scaling}\label{sec:strong}

We next keep the mesh and frequency fixed and increase the GPU
allocation to quantify the reduction in solution time. We use $480^3$
meshes on 1 to 8 GPUs within one node, $960^3$ meshes on 8 to 64 GPUs
across 1 to 8 nodes, and $1920^3$ meshes on 64 to 256 GPUs across 8 to
32 nodes. The $960^3$ and $1920^3$ systems contain approximately 2.66
and 21.3 billion complex unknowns and, in the homogeneous case, span 156
and 316 wavelengths per coordinate direction. The strong-scaling
efficiency $\eta_{\rm s}$ is the speedup in complete solution time
divided by the increase in GPU count, relative to the smallest
allocation for each mesh. Tables~\ref{tab:iofd-strong-480}--\ref{tab:iofd-strong-1920}
list the results.

\begin{table}[htbp]
\centering
\caption{Strong scaling on a fixed $480^3$ mesh using GPUs. Efficiency is based on complete solve time; memory is the maximum measured device usage per GPU.}
\label{tab:iofd-strong-480}
\footnotesize
\begin{adjustbox}{max width=\textwidth}
\begin{tabular}{cccccccc}
\toprule
medium & PC calls & GPUs & $T_{\rm setup}$ & $T_{\rm solve}$ & solve/PC & $\eta_{\rm s}$ & peak/GPU \\
\midrule
homogeneous & 14 & 1 & 0.240~s & 15.895~s & 1.1353~s & 100.0\% & 50.35~GiB \\
homogeneous & 14 & 2 & 0.201~s & 7.888~s & 0.5634~s & 100.8\% & 33.98~GiB \\
homogeneous & 14 & 4 & 0.156~s & 4.442~s & 0.3173~s & 89.5\% & 18.19~GiB \\
homogeneous & 14 & 8 & 0.154~s & 2.617~s & 0.1869~s & 75.9\% & 10.20~GiB \\
\addlinespace
lens & 22 & 1 & 2.335~s & 52.702~s & 2.3955~s & 100.0\% & 48.58~GiB \\
lens & 22 & 2 & 0.239~s & 25.644~s & 1.1656~s & 102.8\% & 36.93~GiB \\
lens & 22 & 4 & 0.177~s & 13.391~s & 0.6087~s & 98.4\% & 19.69~GiB \\
lens & 22 & 8 & 0.158~s & 7.414~s & 0.3370~s & 88.9\% & 10.97~GiB \\
\addlinespace
periodic & 26 & 1 & 2.330~s & 62.929~s & 2.4203~s & 100.0\% & 48.58~GiB \\
periodic & 26 & 2 & 0.252~s & 30.604~s & 1.1771~s & 102.8\% & 36.93~GiB \\
periodic & 26 & 4 & 0.177~s & 15.981~s & 0.6147~s & 98.4\% & 19.69~GiB \\
periodic & 26 & 8 & 0.156~s & 8.834~s & 0.3398~s & 89.0\% & 10.97~GiB \\
\addlinespace
random & 22 & 1 & 2.372~s & 53.242~s & 2.4201~s & 100.0\% & 48.58~GiB \\
random & 22 & 2 & 0.269~s & 25.900~s & 1.1773~s & 102.8\% & 36.93~GiB \\
random & 22 & 4 & 0.186~s & 13.518~s & 0.6144~s & 98.5\% & 19.69~GiB \\
random & 22 & 8 & 0.167~s & 7.474~s & 0.3397~s & 89.0\% & 10.97~GiB \\
\bottomrule
\end{tabular}%
\end{adjustbox}
\end{table}

\begin{table}[htbp]
\centering
\caption{Strong scaling on a fixed $960^3$ mesh using GPUs. Efficiency is based on complete solve time; memory is the maximum measured device usage per GPU.}
\label{tab:iofd-strong-960}
\footnotesize
\begin{adjustbox}{max width=\textwidth}
\begin{tabular}{cccccccc}
\toprule
medium & PC calls & GPUs & $T_{\rm setup}$ & $T_{\rm solve}$ & solve/PC & $\eta_{\rm s}$ & peak/GPU \\
\midrule
homogeneous & 23 & 8 & 0.327~s & 24.706~s & 1.0742~s & 100.0\% & 51.91~GiB \\
homogeneous & 23 & 16 & 0.379~s & 14.455~s & 0.6285~s & 85.5\% & 36.19~GiB \\
homogeneous & 23 & 32 & 0.423~s & 7.839~s & 0.3408~s & 78.8\% & 19.24~GiB \\
homogeneous & 23 & 64 & 0.309~s & 4.544~s & 0.1976~s & 68.0\% & 10.73~GiB \\
\addlinespace
lens & 45 & 8 & 2.673~s & 98.500~s & 2.1889~s & 100.0\% & 50.15~GiB \\
lens & 45 & 16 & 0.397~s & 51.275~s & 1.1395~s & 96.1\% & 37.95~GiB \\
lens & 46 & 32 & 0.406~s & 28.234~s & 0.6138~s & 87.2\% & 19.99~GiB \\
lens & 45 & 64 & 0.288~s & 15.964~s & 0.3548~s & 77.1\% & 10.99~GiB \\
\addlinespace
periodic & 50 & 8 & 3.181~s & 110.465~s & 2.2093~s & 100.0\% & 50.15~GiB \\
periodic & 50 & 16 & 0.401~s & 57.482~s & 1.1496~s & 96.1\% & 37.95~GiB \\
periodic & 50 & 32 & 0.412~s & 30.895~s & 0.6179~s & 89.4\% & 19.99~GiB \\
periodic & 50 & 64 & 0.297~s & 17.933~s & 0.3587~s & 77.0\% & 10.99~GiB \\
\addlinespace
random & 41 & 8 & 3.747~s & 90.591~s & 2.2095~s & 100.0\% & 50.15~GiB \\
random & 41 & 16 & 0.415~s & 47.188~s & 1.1509~s & 96.0\% & 37.95~GiB \\
random & 41 & 32 & 0.419~s & 25.437~s & 0.6204~s & 89.0\% & 19.99~GiB \\
random & 41 & 64 & 0.301~s & 14.608~s & 0.3563~s & 77.5\% & 10.99~GiB \\
\bottomrule
\end{tabular}%
\end{adjustbox}
\end{table}

\begin{table}[htbp]
\centering
\caption{Strong scaling on a fixed $1920^3$ mesh using GPUs. Efficiency is based on complete solve time; memory is the maximum measured device usage per GPU.}
\label{tab:iofd-strong-1920}
\footnotesize
\begin{adjustbox}{max width=\textwidth}
\begin{tabular}{cccccccc}
\toprule
medium & PC calls & GPUs & $T_{\rm setup}$ & $T_{\rm solve}$ & solve/PC & $\eta_{\rm s}$ & peak/GPU \\
\midrule
homogeneous & 42 & 64 & 4.922~s & 49.375~s & 1.1756~s & 100.0\% & 46.83~GiB \\
homogeneous & 42 & 128 & 0.454~s & 26.273~s & 0.6256~s & 94.0\% & 36.69~GiB \\
homogeneous & 42 & 256 & 0.339~s & 16.638~s & 0.3961~s & 74.2\% & 19.04~GiB \\
\addlinespace
lens & 106 & 64 & 4.708~s & 239.432~s & 2.2588~s & 100.0\% & 50.25~GiB \\
lens & 102 & 128 & 0.468~s & 116.840~s & 1.1455~s & 102.5\% & 38.20~GiB \\
lens & 104 & 256 & 0.368~s & 69.695~s & 0.6701~s & 85.9\% & 20.61~GiB \\
\addlinespace
periodic & 113 & 64 & 4.455~s & 257.479~s & 2.2786~s & 100.0\% & 50.25~GiB \\
periodic & 114 & 128 & 0.455~s & 131.711~s & 1.1554~s & 97.7\% & 38.20~GiB \\
periodic & 113 & 256 & 0.371~s & 76.439~s & 0.6765~s & 84.2\% & 20.61~GiB \\
\addlinespace
random & 84 & 64 & 4.382~s & 191.471~s & 2.2794~s & 100.0\% & 50.25~GiB \\
random & 84 & 128 & 0.533~s & 97.003~s & 1.1548~s & 98.7\% & 38.20~GiB \\
random & 85 & 256 & 0.381~s & 57.327~s & 0.6744~s & 83.5\% & 20.61~GiB \\
\bottomrule
\end{tabular}%
\end{adjustbox}
\end{table}

The strong-scaling results reduce billion-unknown Maxwell solves to
seconds. On the $480^3$ mesh, 8 GPUs complete every solve in at most
8.8 seconds. On the $960^3$ mesh, increasing the allocation from 8 to 64
GPUs shortens the homogeneous solve from 24.7 to 4.5 seconds and the
heterogeneous solves from 90.6--110.5 to 14.6--17.9 seconds. On the
$1920^3$ mesh, increasing the allocation from 64 to 256 GPUs shortens
the homogeneous solve, which spans 316 wavelengths per direction, from
49.4 to 16.6 seconds and the heterogeneous solves from 191.5--257.5 to
57.3--76.4 seconds, with strong-scaling efficiencies of 74.2--85.9\%.
\FloatBarrier

\subsection{CPU scaling}\label{sec:cpu}

The CPU implementation carries the same discretization and three-grid method to
full-wave Maxwell computation beyond 500 billion unknowns. The largest
homogeneous problem is solved on a $5760^3$ mesh with approximately
574 billion complex unknowns, and its domain spans 956 wavelengths in
each coordinate direction. The Green-field comparison in
Table~\ref{tab:iofd-green} establishes the accuracy of this solution;
here we examine the resources and solution time needed to reach it.

The CPU implementation runs eight MPI ranks per two-socket Intel Xeon
Platinum 8358P node, each rank using eight OpenMP threads, so that the
64 cores of a node are filled; the ranks form a three-dimensional brick
partition. CPU memory denotes the total allocation over all ranks. The
weak-scaling sequence pairs $960^3$, $1920^3$, $3840^3$, and $5760^3$
meshes with 1, 8, 64, and 216 nodes, so that every rank holds a $480^3$
block. We use the per-application efficiency
\eqref{eq:iofd-weak-efficiency} with the one-node $960^3$ case as
reference; the strong-scaling efficiency is based on complete solution
time relative to the smallest allocation for each mesh.

The per-application weak scaling remains nearly ideal over the 216-fold
increase in mesh volume and node count (Table~\ref{tab:iofd-cpu-weak}).
The application cost rises only from 20.4 to 21.4 seconds, giving
95.7\% efficiency on 13\,824 cores.

This sustained throughput turns the 574-billion-unknown system into a
tractable complete solve. The problem spanning 956 wavelengths per
direction converges in 119 applications and approximately 42 minutes,
after a hierarchy setup of about 6 seconds. Because the storage grows
in proportion to the mesh volume, memory per node stays at about
365~GiB from 1 to 216 nodes, so larger electrical domains are reached by
adding nodes at a fixed storage budget per node and a nearly fixed cost
per application.

\begin{table}[htbp]
\centering
\caption{CPU weak scaling for the homogeneous medium, relative to the one-node $960^3$ case.}
\label{tab:iofd-cpu-weak}
\footnotesize
\begin{adjustbox}{max width=\textwidth}
\begin{tabular}{cccccccccc}
\toprule
mesh & unknowns & nodes & CPU cores & PC calls & setup & solve & solve/PC & $\eta_{\rm w}$ & memory \\
\midrule
$960^3$ & $2.66\times10^{9}$ & 1 & 64 & 23 & 5.32~s & 470.21~s & 20.444~s & 100.0\% & 364.6~GiB \\
$1920^3$ & $2.13\times10^{10}$ & 8 & 512 & 42 & 6.04~s & 874.55~s & 20.823~s & 98.2\% & 2918.1~GiB \\
$3840^3$ & $1.70\times10^{11}$ & 64 & 4096 & 80 & 6.62~s & 1658.30~s & 20.729~s & 98.6\% & 23350.5~GiB \\
$5760^3$ & $5.74\times10^{11}$ & 216 & 13824 & 119 & 6.09~s & 2541.28~s & 21.355~s & 95.7\% & 78837.3~GiB \\
\bottomrule
\end{tabular}%
\end{adjustbox}
\end{table}

To complement this increase in problem size, we fix the mesh and add
nodes to shorten the complete solve (Table~\ref{tab:iofd-cpu-strong}). For $3840^3$, a system with
approximately 170 billion complex unknowns, increasing the allocation
from 64 to 256 nodes reduces the solve from 1658 to 447 seconds, with an
efficiency of 92.8\%. An eightfold increase in node count reduces the
$1920^3$ solve from 875 to 124 seconds and the $960^3$ solve from 470 to
68 seconds, with efficiencies of 88.1\% and 87.1\%, respectively.

\begin{table}[htbp]
\centering
\caption{CPU strong scaling for the homogeneous medium, relative to the smallest allocation for each mesh.}
\label{tab:iofd-cpu-strong}
\footnotesize
\begin{adjustbox}{max width=\textwidth}
\begin{tabular}{cccccccccc}
\toprule
mesh & unknowns & nodes & CPU cores & PC calls & setup & solve & solve/PC & $\eta_{\rm s}$ & memory \\
\midrule
$960^3$ & $2.66\times10^{9}$ & 1 & 64 & 23 & 5.32~s & 470.21~s & 20.444~s & 100.0\% & 364.6~GiB \\
$960^3$ & $2.66\times10^{9}$ & 2 & 128 & 23 & 2.74~s & 241.09~s & 10.482~s & 97.5\% & 370.0~GiB \\
$960^3$ & $2.66\times10^{9}$ & 4 & 256 & 23 & 1.60~s & 124.47~s & 5.412~s & 94.4\% & 376.0~GiB \\
$960^3$ & $2.66\times10^{9}$ & 8 & 512 & 23 & 0.82~s & 67.50~s & 2.935~s & 87.1\% & 386.9~GiB \\
\addlinespace
$1920^3$ & $2.13\times10^{10}$ & 8 & 512 & 42 & 6.04~s & 874.55~s & 20.823~s & 100.0\% & 2918.1~GiB \\
$1920^3$ & $2.13\times10^{10}$ & 16 & 1024 & 42 & 3.07~s & 441.70~s & 10.517~s & 99.0\% & 2959.5~GiB \\
$1920^3$ & $2.13\times10^{10}$ & 32 & 2048 & 42 & 1.63~s & 227.66~s & 5.421~s & 96.0\% & 3011.1~GiB \\
$1920^3$ & $2.13\times10^{10}$ & 64 & 4096 & 42 & 0.83~s & 124.14~s & 2.956~s & 88.1\% & 3097.1~GiB \\
\addlinespace
$3840^3$ & $1.70\times10^{11}$ & 64 & 4096 & 80 & 6.62~s & 1658.30~s & 20.729~s & 100.0\% & 23350.5~GiB \\
$3840^3$ & $1.70\times10^{11}$ & 128 & 8192 & 80 & 3.36~s & 853.87~s & 10.673~s & 97.1\% & 23685.4~GiB \\
$3840^3$ & $1.70\times10^{11}$ & 256 & 16384 & 80 & 1.65~s & 446.75~s & 5.584~s & 92.8\% & 24117.7~GiB \\
\bottomrule
\end{tabular}%
\end{adjustbox}
\end{table}

\FloatBarrier

\section{Conclusions}\label{sec:conclusions}

We have developed a Maxwell IOFD discretization and a three-grid solver
that together make accurate three-dimensional full-wave computation
practical over hundreds of wavelengths. Compact face and edge filters
retain the Yee curl structure while providing enough freedom to control
both transverse dispersion branches. The symbol analysis connects the
filter coefficients to the numerical wavenumbers and bounds the effects
of coefficient interpolation. Optimizing at the operating resolution
then suppresses the phase error that would otherwise accumulate along
long propagation paths. This improvement is realized in the computed
fields: at only 6 points per wavelength, the homogeneous problem remains
accurate on a domain 956 wavelengths wide in every coordinate direction,
with a vector $L^2$ error of $4.38\%$.

The three-grid hierarchy converts this low-dispersion representation into
an efficient large-scale solver. An unshifted intermediate-grid equation
supplies the propagating correction, while a complex-shifted auxiliary
cycle makes its iterative solution inexpensive. The resulting compact,
matrix-free operations keep storage proportional to the number of
unknowns and support effective parallel execution. Complete solves with
approximately 71.7 billion complex unknowns take 77--504 seconds on
216 GPUs across homogeneous, focusing, periodic, and random
media. On CPUs, the same method reaches approximately
574 billion unknowns with per-application weak-scaling efficiencies above
$95\%$. These results establish the value of combining dispersion
optimization with a fully matrix-free multilevel solver: accurate
propagation at modest resolution can be carried through to fast solutions
of very large three-dimensional Maxwell systems.

\appendix
\section{Dispersion Analysis}\label{sec:dispersion}

In this appendix, we derive the dispersion properties behind the
coefficient design of Section~\ref{sec:symbol}: the leading phase errors
of the two polarizations (Appendix~\ref{sec:leading}), the radial control
of the wavenumber error under coefficient perturbation and interpolation
(Appendix~\ref{sec:radial}), and the residues that the source and field
filter corrects (Appendix~\ref{sec:residues}). In the notation of \eqref{eq:dimensionlesssymbol},
the stiffness symbol is
\begin{equation}
 \widehat K=\begin{pmatrix}
 f_2s_3^2+f_3s_2^2&-f_3s_1s_2&-f_2s_1s_3\\
 -f_3s_1s_2&f_1s_3^2+f_3s_1^2&-f_1s_2s_3\\
 -f_2s_1s_3&-f_1s_2s_3&f_1s_2^2+f_2s_1^2
 \end{pmatrix}.
 \label{eq:Kentries}
\end{equation}
For a prescribed frequency, $\eta$ fixes the coefficients $\cc(\eta)$
while $\tht$ ranges over the Fourier modes. The two roots in
\eqref{eq:metrics} also control the artificial difference between the
polarization wavenumbers,
\begin{equation}
 \frac{|k_{h,1}-k_{h,2}|}{k}
     =|q_1-q_2|\le|\delta_1|+|\delta_2|
     \le2\mathcal E_{\ph}(\eta).
 \label{eq:polarizationdifference}
\end{equation}
and cubic symmetry allows the angular search to be restricted to a
symmetry sector.

\subsection{Leading phase errors}\label{sec:leading}
\begin{proposition}[Leading phase errors of the two polarizations]
\label{prop:expansion}
Fix a real coefficient vector $\cc$ and a unit direction $\dvec$.
Let $U\in\R^{3\times2}$ have orthonormal columns spanning
$\dvec^\perp$, and define
\begin{equation}
 S_4=\sum_{j=1}^3d_j^4,\qquad
 D_T=U^T\diag(d_1^2,d_2^2,d_3^2)U.
 \label{eq:transverseD}
\end{equation}

The two positive generalized eigenvalues of
$(\widehat K(t\dvec;\cc),\widehat M(t\dvec;\cc))$ satisfy, as $t\to0$,
\begin{equation}
 \lambda_b(t\dvec;\cc)=t^2+t^4\gamma_b+O(t^6),
 \label{eq:asymptoticlambda}
\end{equation}
where $\gamma_1,\gamma_2$ are the eigenvalues of
\begin{equation}
 \Gamma_T=\left[\frac{b_t-a}{2}
            +\left(\frac a2-\frac1{12}\right)S_4\right]I_2
            +\frac{a+b_n-b_t}{2}D_T.
 \label{eq:Gamma}
\end{equation}
The near-physical roots therefore obey
\begin{equation}
 \delta_b=-\frac{\kappa^2}{2}\gamma_b+O(\kappa^4).
 \label{eq:phaseexpansion}
\end{equation}
The leading polarization difference vanishes for every direction when
$a+b_n-b_t=0$. The $O(\kappa^2)$ phase-error term vanishes for both
polarizations in every direction if and only if
\begin{equation}
 a=\frac16,\qquad b_n=0,\qquad b_t=\frac16.
 \label{eq:leadingcancel}
\end{equation}
\end{proposition}
Applying this expansion to the Yee member $\cc=0$ gives
\begin{equation}
 \lambda_1=\lambda_2=\sum_j4\sin^2(\theta_j/2),\qquad
 \delta_b=\frac{\kappa^2}{24}\sum_jd_j^4+O(\kappa^4).
 \label{eq:yeephase}
\end{equation}
which is the leading, direction-dependent Yee error removed by
\eqref{eq:leadingcancel}. The remaining parameters $b_{tt}$ and $b_{nt}$
first enter at order $t^6$, and $b_{ntt}$ at order $t^8$.

\subsection{Radial error control and interpolation}\label{sec:radial}
Admissibility also bounds the spread of the constitutive symbols,
\begin{equation}
 \chi(\cc)=\max\left\{
 \frac{\max_{r,\tht}f_r}{\min_{r,\tht}f_r},
 \frac{\max_{r,\tht}m_r}{\min_{r,\tht}m_r}\right\},
 \label{eq:spread}
\end{equation}
with all extrema attained at vertices. Evaluating a branch on the
physical sphere gives the residual
\begin{equation}
 r_b(\dvec;\eta,\cc)
       =\frac{\lambda_b(\kappa\dvec;\cc)}{\kappa^2}-1.
 \label{eq:sphereresidual}
\end{equation}
whose effect on the wavenumber depends on how rapidly the branch crosses
the level $\kappa^2$.

\begin{proposition}[Radial phase-error control]\label{prop:radial}
Fix $\kappa>0$, a unit direction $\dvec$, and an admissible coefficient
vector $\cc$. Let $\lambda_b$ be either positive generalized eigenvalue
of $(\widehat K,\widehat M)$, and define
$F_b(q)=\lambda_b(q\kappa\dvec;\cc)-\kappa^2$ on
$I=[q_-,q_+]$. Suppose that $0<q_-<1<q_+$, that
$q\kappa\dvec\in[-\pi,\pi]^3$ for every $q\in I$, and that
$F_b(q_-)\le0\le F_b(q_+)$. If
\begin{equation}
 F_b(q_2)-F_b(q_1)\ge\gamma_{\rm r}\kappa^2(q_2-q_1),
 \qquad q_2>q_1,\quad q_1,q_2\in I,
 \label{eq:radialslope}
\end{equation}
for some $\gamma_{\rm r}>0$, then $F_b$ has a unique root $q_b\in I$, and
\begin{equation}
 |q_b-1|\le\frac{|r_b|}{\gamma_{\rm r}}.
 \label{eq:rootbound}
\end{equation}
At a simple generalized eigenpair normalized by
$\vct p_b^*\widehat M\vct p_b=1$, the radial derivative is
\begin{equation}
 F_b'(q)=\vct p_b^*
    (\partial_q\widehat K-\lambda_b\partial_q\widehat M)\vct p_b,
 \qquad\partial_q=\kappa\dvec\cdot\nabla_{\tht}.
 \label{eq:HF}
\end{equation}
\end{proposition}
A differentiable branch with $F_b'(q)\ge\gamma_{\rm r}\kappa^2$ satisfies
\eqref{eq:radialslope}, and linearizing at $q=1$ gives
\begin{equation}
 \delta_b\simeq
 -\frac{\lambda_b(\kappa\dvec)-\kappa^2}
        {\kappa\dvec\cdot\nabla_{\tht}\lambda_b(\kappa\dvec)}.
 \label{eq:newtonestimate}
\end{equation}

\begin{proposition}[Coefficient perturbations and wavenumber error]
\label{prop:perturbation}
Fix $\kappa>0$ and a unit direction $\dvec$. Let $\cc$ and
$\widetilde\cc$ be admissible coefficient vectors, and let
$\lambda_b,\widetilde\lambda_b$ denote their ordered positive generalized
eigenvalues. On a common interval $I$ of radial arguments
$\tht=q\kappa\dvec$ in the Fourier cube, suppose that
$\widehat M,\widetilde{\widehat M}\succeq\underline m I_3$, with
$\underline m>0$, and both stiffness symbols have norm at most
$K_{\max}$ in the spectral norm. Put $\Delta\cc=\widetilde\cc-\cc$ and
\begin{equation}
 L_M=2|\Delta b_n|+4|\Delta b_t|+4|\Delta b_{tt}|
              +8|\Delta b_{nt}|+8|\Delta b_{ntt}|.
 \label{eq:LM}
\end{equation}
At each such Fourier argument, the ordered eigenvalues satisfy
\begin{equation}
 |\widetilde\lambda_b-\lambda_b|
 \le\frac{2|\svec|^2}{\underline m}|\Delta a|
       +\frac{K_{\max}}{\underline m^2}L_M
       =:\varepsilon_\lambda(\tht).
 \label{eq:eigenperturb}
\end{equation}
If the corresponding roots $q_b,\widetilde q_b$ lie in that interval
and the unperturbed branch satisfies \eqref{eq:radialslope}, then
\begin{equation}
 |\widetilde q_b-q_b|
       \le\frac{\sup_{q\in I}\varepsilon_\lambda(q\kappa\dvec)}
                    {\gamma_{\rm r}\kappa^2}.
 \label{eq:phaseperturb}
\end{equation}
\end{proposition}
For a twice continuously differentiable reference coefficient curve
$\cc^\star$, linear interpolation has the componentwise error bound
\begin{equation}
 \|c_i^\star-\mathcal I c_i^\star\|_{L^\infty[\eta_j,\eta_{j+1}]}
 \le\frac{(\eta_{j+1}-\eta_j)^2}{8}
            \|(c_i^\star)''\|_{L^\infty[\eta_j,\eta_{j+1}]}.
 \label{eq:linearinterpolationerror}
\end{equation}
Under the conditions of Proposition~\ref{prop:perturbation}, substitution
into \eqref{eq:LM}--\eqref{eq:phaseperturb} gives a second-order
interpolation contribution to the wavenumber error when the mass lower
bound and radial slope remain bounded away from zero.

\subsection{Residues and the source and field filter}\label{sec:residues}
For admissible coefficients let $\vct p_0,\vct p_1,\vct p_2$ be
$\widehat M$-orthonormal generalized eigenvectors of
$(\widehat K,\widehat M)$, with $\lambda_0=0$. Then
\begin{equation}
 \mathscr P^{-1}
 =-\frac{\vct p_0\vct p_0^*}{\kappa^2}
       +\sum_{b=1}^2\frac{\vct p_b\vct p_b^*}{\lambda_b-\kappa^2}.
 \label{eq:discreteresolvent}
\end{equation}

\begin{proposition}[The residue of a simple transverse pole]
\label{prop:residue}
Fix $\eta$ and an admissible coefficient vector. Suppose that the matrix
symbol $\mathscr P(\tht_*;\eta)$ in \eqref{eq:dimensionlesssymbol} has a simple zero
eigenvalue with Euclidean unit eigenvector $e_b$. Choose a unit normal
$n_b$ to its characteristic surface and assume
\begin{equation}
 \alpha_b=e_b^*(n_b\cdot\nabla_{\tht}\mathscr P)e_b\ne0
       \quad\text{at }\tht_*.
 \label{eq:alpha}
\end{equation}
Then
\begin{equation}
 \mathscr P(\tht_*+tn_b;\eta)^{-1}
       =\frac{e_be_b^*}{t\alpha_b}+O(1),\qquad t\to0.
 \label{eq:pole}
\end{equation}
\end{proposition}
In the generalized normalization of \eqref{eq:discreteresolvent},
$e_b=\vct p_b/\|\vct p_b\|_2$ and
\begin{equation}
 \alpha_b
 =\frac{n_b\cdot\nabla_{\tht}\lambda_b}{\|\vct p_b\|_2^2}.
 \label{eq:ordinarygradient}
\end{equation}
on the characteristic surface. The continuous transverse poles have
strength $1/(2\kappa)$, so multiplying the inverse symbol by
$q_Q^2\approx\alpha_b/(2\kappa)$ moves the discrete residues toward the
continuous ones. The filter \eqref{eq:Qweights} has the Fourier symbol
\begin{equation}
 \begin{split}
 q_Q(\tht)={}&1+\beta_1\left(\sum_j\cos\theta_j-3\right)
       +\beta_2\left(\sum_{i<j}\cos\theta_i\cos\theta_j-3\right)\\
       &+\beta_3\left(\prod_j\cos\theta_j-1\right),
 \end{split}
 \label{eq:Qsymbol}
\end{equation}
with $q_Q(0)=1$.

\section{Proofs for Section~\ref{sec:discretization}}\label{sec:proofs}
\subsection*{Proof of Proposition~\ref{prop:consistency}}
\begin{proof}
For a smooth component field $w$, Taylor expansion on its lattice gives
$\Delta_jw=\tfrac12h^2\partial_j^2w+O(h^4)$. Thus
$F_h=I+O(h^2)$ and $M_h=I+O(h^2)$ on smooth fields. The staggered
differences in $C_h$ and its reverse curl are second-order accurate at
their respective locations, proving consistency.

The two stiffness identities follow from $C_hG_h=0$ and its adjoint.
Applying $G_h^*$ to
the discrete equation proves \eqref{eq:discretegauss}.
\end{proof}

\subsection*{Proof of Proposition~\ref{prop:branches}}
\begin{proof}
For nonzero $\svec$, the cross-product matrix $S$ has rank two and kernel
$\operatorname{span}\{\svec\}$. Since $\widehat F$ is positive definite,
$\widehat K=S^T\widehat F S$ has the same kernel and two positive
eigenvalues after congruence by $\widehat M^{-1/2}$. Multiplying
$\widehat K\vct p=\lambda\widehat M\vct p$ by $\svec^T$ gives
the transverse constraint for $\lambda>0$.

The two nonzero eigenvalues
of $B$ have sum $\tr B$ and squared sum $\tr(B^2)$, which proves
\eqref{eq:explicitbranches}. Factoring
$\mathscr P=\widehat M^{1/2}(B-\kappa^2I)\widehat M^{1/2}$ gives
\eqref{eq:characteristicfactorization}.

Finally, $f_r$ is affine in $z_r$,
and $m_r$ is affine in each $z_j$ separately. Each value inside the cube
is therefore a convex combination of vertex values.
\end{proof}

\subsection*{Proof of Proposition~\ref{prop:expansion}}
\begin{proof}
Let $S_d=[\dvec]_\times$ and
$S_{d^3}=[(d_1^3,d_2^3,d_3^3)^T]_\times$. Expanding sine and cosine gives
\begin{align*}
 S(t\dvec)&=tS_d-\frac{t^3}{24}S_{d^3}+O(t^5),\\
 \widehat F(t\dvec)&=I-\frac{at^2}{2}\diag(d_j^2)+O(t^4),\\
 \widehat M(t\dvec)&=I+t^2M_2+O(t^4),\\
 M_2&=-\frac12\bigl[b_tI+(b_n-b_t)\diag(d_j^2)\bigr].
\end{align*}
Consequently,
$\widehat K=t^2(I-\dvec\dvec^T)+t^4K_4+O(t^6)$, where
\begin{equation}
 K_4=-\frac1{24}(S_{d^3}^TS_d+S_d^TS_{d^3})
                  -\frac a2 S_d^T\diag(d_j^2)S_d.
 \label{eq:K4}
\end{equation}
To extract the transverse correction, expand
$B=\widehat M^{-1/2}\widehat K\widehat M^{-1/2}$.
Its leading eigenspace is $\dvec^\perp$, and restriction to that space gives
the $t^4$ correction $U^T(K_4-M_2)U$. The transverse identities are
\begin{align*}
 U^T(S_{d^3}^TS_d+S_d^TS_{d^3})U&=2S_4I_2,\\
 U^TS_d^T\diag(d_j^2)S_dU&=(1-S_4)I_2-D_T.
\end{align*}
The second identity follows from the action of $S_d$ as a right-angle
rotation on $\dvec^\perp$. Substitution gives \eqref{eq:Gamma}.

Inserting $q_b=1+O(\kappa^2)$ into
$\lambda_b(q_b\kappa\dvec)=\kappa^2$ gives
\eqref{eq:phaseexpansion}. The scalar term of $\Gamma_T$ is common to both
polarizations; its $D_T$ term vanishes when $a+b_n-b_t=0$.

For all-direction cancellation of the leading phase error,
$\dvec=(1,0,0)^T$ forces $b_t=1/6$. Taking
$\dvec=(1,1,0)^T/\sqrt2$, for which $D_T$ has eigenvalues $0$ and $1/2$,
then forces $a=1/6$ and $b_n=0$. These values make $\Gamma_T$ zero
for every direction, proving sufficiency.
\end{proof}

\subsection*{Proof of Proposition~\ref{prop:radial}}
\begin{proof}
The generalized eigenvalues are continuous because the mass symbol is
positive definite. The sign condition and intermediate-value theorem give
a root; \eqref{eq:radialslope} makes it unique. Applying that inequality
between $1$ and $q_b$ gives
$\gamma_{\rm r}\kappa^2|q_b-1|\le|F_b(1)|=\kappa^2|r_b|$, which proves
\eqref{eq:rootbound}.

For the derivative formula, differentiate
$\widehat K\vct p_b=\lambda_b\widehat M\vct p_b$ and multiply by
$\vct p_b^*$. The eigenvector-derivative terms cancel, yielding
\eqref{eq:HF}.
\end{proof}

\subsection*{Proof of Proposition~\ref{prop:perturbation}}
\begin{proof}
The bounds $|z_j|\le2$ give
$\|\widetilde{\widehat M}-\widehat M\|_2\le L_M$.
Since $\|S\|_2=|\svec|$,
$\|\widetilde{\widehat K}-\widehat K\|_2
 \le2|\svec|^2|\Delta a|$.
For any Euclidean unit vector $v$, compare the generalized Rayleigh
quotients with numerators $v^*\widehat Kv$ and
$v^*\widetilde{\widehat K}v$. Their difference is bounded by
$\|\Delta\widehat K\|_2/\underline m
 +K_{\max}\|\Delta\widehat M\|_2/\underline m^2$.
The min--max principle gives \eqref{eq:eigenperturb} for the ordered
eigenvalues, including at a polarization coincidence.

At
$\widetilde q_b$, the perturbed eigenvalue is $\kappa^2$.
Applying \eqref{eq:radialslope} between the two roots and using
\eqref{eq:eigenperturb} proves \eqref{eq:phaseperturb}.
\end{proof}

\subsection*{Proof of Proposition~\ref{prop:residue}}
\begin{proof}
The simple eigenvalue of the Hermitian matrix $\mathscr P$ is
$t\alpha_b+O(t^2)$ along this crossing. Its spectral projector is
$e_be_b^*+O(t)$, while the other eigenvalues remain separated from zero.
Inverting the spectral decomposition proves \eqref{eq:pole}.
\end{proof}

\section{The Random Material Field}\label{app:random}

The fixed realization in \eqref{eq:randommodel} is the trigonometric sum
\begin{equation}
 R(\x)=\sum_{j=1}^8 A_j\,\psi_j(2\pi\vct\nu_j\cdot\x+\varphi_j),
 \label{eq:randomfourier}
\end{equation}
where each $\psi_j$ is either sine or cosine, as specified together with
the amplitudes, frequencies, and phases in Table~\ref{tab:randomfield}.
This representation defines the same material at every grid resolution.

\begin{table}[htbp]
\centering\small
\caption{Fourier coefficients of the random material field.}
\label{tab:randomfield}
\begin{tabular}{cccc}
\toprule
$A_j$ & $\vct\nu_j$ & $\varphi_j$ & $\psi_j$\\
\midrule
0.71 & $(0.50,1.00,-0.25)$ & 0.37 & $\cos$\\
0.53 & $(1.25,-0.50,0.75)$ & 1.11 & $\sin$\\
0.47 & $(-0.75,1.50,0.50)$ & 2.23 & $\cos$\\
0.41 & $(1.50,0.25,-1.00)$ & 3.17 & $\sin$\\
0.37 & $(0.25,-1.25,1.50)$ & 4.03 & $\cos$\\
0.31 & $(-1.50,0.75,1.25)$ & 5.19 & $\sin$\\
0.29 & $(1.00,1.00,1.00)$ & 0.83 & $\cos$\\
0.23 & $(-1.25,-0.75,0.50)$ & 2.71 & $\sin$\\
\bottomrule
\end{tabular}
\end{table}
\FloatBarrier

\section*{Declaration of AI assistance}
AI tools were used to assist with the proofs, the
implementation of the algorithms, and the writing of the manuscript. The
authors take full responsibility for the content of this paper.

\end{document}